# Fourier-Stieltjes Algebras of Locally Compact Groupoids

by

Arlan Ramsay
and
Martin E. Walter

ABSTRACT. For locally compact groups, Fourier algebras and Fourier-Stieltjes algebras have proved to be useful dual objects. They encode the representation theory of the group via the positive definite functions on the group: positive definite functions correspond to cyclic representations and span these algebras as linear spaces. They encode information about the algebra of the group in the geometry of the Banach space structure, and the group appears as a topological subspace of the maximal ideal space of the algebra [W1, W2]. Because groupoids and their representations appear in studying operator algebras, ergodic theory, geometry, and the representation theory of groups, it would be useful to have a duality theory for them. This paper gives a first step toward extending the theory of Fourier-Stieltjes algebras from groups to groupoids. If $G$ is a locally compact (second countable) groupoid, we show that $\mathcal{B}(G)$, the linear span of the Borel positive definite functions on $G$, is a Banach algebra when represented as an algebra of completely bounded maps on a $C^*$-algebra associated with $G$. This necesarily involves identifying equivalent elements of $\mathcal{B}(G)$. An example shows that the linear span of the continuous positive definite functions need not be complete. For groups, $\mathcal{B}(G)$ is isometric to the Banach space dual of $C^*(G)$. For groupoids, the best analog of that fact is to be found in a representation of $\mathcal{B}(G)$ as a Banach space of completely bounded maps from a $C^*$-algebra associated with $G$ to a $C^*$-algebra associated with the equivalence relation induced by $G$. This paper adds weight to the clues in the earlier study of Fourier-Stieltjes algebras that there is a much more general kind of duality for Banach algebras waiting to be explored [W5].

## Introduction

As suggested by the title, this paper connects two lines of earlier work, and we begin with an abbreviated history of each of these lines, in order of appearance. After the history, we will state our main results and outline the body of the paper. We mention here that some basic definitions can be found in Section 1 and that we assume locally compact spaces are second countable. More background on groupoids is available in [C1, C2, M1, M2, H1, H2, R1, R3, Re1, Re2, S1, S2]. The necessary background on Fourier-Stieltjes algebras can be obtained from [E, W1, W2].

In 1963, G.W. Mackey, [M1], introduced the notion of virtual group as a tool and context for several kinds of problems in analysis and geometry. Virtual groups are (equivalence classes of) groupoids having suitable measure theoretic structure and the property of ergodicity. Ergodicity makes a groupoid more group-like, but many results on groupoids do not require ergodicity. Among the structures which fit naturally into the study of groupoids are groups, group actions, equivalence relations (including foliations), ordinary spaces, and examples made from these by restricting to a part of the underlying space.

The original motivation for studying groupoids was provided by Mackey's theory of unitary representations of group extensions. The idea has been applied to that subject: for example see [B-M-R, B-C-M-R, R2]. In his original paper, Mackey also showed the relevance of the idea for ergodic group actions in general, and a number of applications have been made there, for example see [M2, M–Z, Z1, Z2, Z3, Z4].

Most uses of groupoids have been in the study of operator algebras, another approach to understanding and exploiting symmetry. Several pioneering papers should be mentioned. P. Hahn proved the existence of Haar measures for measured groupoids, whether ergodic or not, and used this to make convolution algebras and study von Neumann algebras associated with measured groupoids [H1, H2, H3]. J. Feldman and C.C. Moore made a thorough analysis of ergodic equivalence relations that have countable equivalence classes, showing

The first author was partially supported by National Science Foundation grants MCS-8241056 and DMS-8404104 and by a University of Colorado Faculty Fellowship. Both authors did some work on the paper at the Mathematical Sciences Research Institute, Berkeley.

that the von Neumann algebras attached to them are exactly the factors that have Cartan subalgebras [F-M1, F-M2]. A. Connes introduced a variation on the approach of Mackey in [C1], in particular by working without a chosen invariant measure class. This approach has some advantages for applications to foliations and to $C^*$-algebras [C2]. J. Renault [Re1] studied $C^*$-algebras generated by convolution algebras on locally compact groupoids endowed with Haar systems, not using invariant measure classes. It was shown in [R3] that measured groupoids may be assumed to have locally compact topologies, with no loss in generality. (The continuity of the Haar system was neglected in that paper, but can be achieved.) Thus the study of operator algebras associated with groupoid symmetry can always be confined to locally compact groupoids, whether one is interested in $C^*$-algebras or von Neumann algebras.

Basically one can say that locally compact groupoids occur in situations where there is symmetry that is made evident by the presence of an equivalence relation. Many of these are associated either with group actions or foliations. It can be surprising how group-like both group actions and foliations can be. In particular, some of the papers mentioned above have included information about the unitary representations of groupoids. However, there is no treatment of duality theory for groupoids, and we intend to make a beginning here.

In 1964 P. Eymard, [E], introduced Fourier and Fourier-Stieltjes algebras for non-commutative locally compact groups. Roughly, the Fourier-Stieltjes algebra of a locally compact group, $G$, denoted $B(G)$, is the unitary representation theory of $G$ equipped with some additional algebraic and geometric structure. More precisely, $B(G)$ is the set of finite linear combinations of continuous positive definite functions on $G$ equipped with a norm, which makes $B(G)$ a *commutative* Banach algebra (using *pointwise* addition and multiplication). The elements of $B(G)$ are exactly the matrix entries of unitary representations of $G$. A primary source of intuition is the fact that when $G$ is abelian, $B(G)$ is the isometric, inverse Fourier-Stieltjes transform of $M^1(\hat{G})$, the convolution, Banach algebra of finite, regular Borel measures on $\hat{G}$, the dual group (of characters) of $G$. Thus $B(G)$, as a Banach algebra, "is" $M^1(\hat{G})$. The fact that $B(G)$ exists (as a *commutative* Banach algebra) when $G$ is not abelian leads one to hope that a useful duality theory exists for non-abelian groups which is in spirit similar to the application rich Pontriagin-Van Kampen duality for abelian locally compact groups. That such a duality theory exists has been established by M.E. Walter in [W2] by proving that:

(1) $B(G)$ is a complete invariant of $G$, i.e., $B(G_1)$ and $B(G_2)$ are isometrically isomorphic as Banach algebras, if and only if $G_1$ and $G_2$ are topologically isomorphic as locally compact groups, and

(2) There is an explicit process for recovering $G$ given its "dual object", $B(G)$. Exactly how useful this theory will be remains to be seen since all but a few of the hoped for important applications await rigorous proof.

For various reasons it turns out that it may be more fruitful to look at $B(G)$ from a broader perspective than that afforded by the category of locally compact groups. Namely, in [W3] it is seen that there is a natural duality theory for a "large" collection of Banach algebras that extends in a precise way the Pontriagin duality for abelian groups as well as the above-mentioned duality for non-abelian groups. The theory of $C^*$-algebras plays a large role both technically and intuitively in this duality theory.

In an effort to understand this new duality theory better, as well as to generate meaningful applications and examples of a concrete nature, in this paper we have answered affirmatively the question: Does a locally compact groupoid $G$ have a Fourier-Stieltjes algebra? For groupoids, there is more than one candidate for the Fourier-Stieltjes algebra, and the details are more technical than for groups, but there is an affirmative answer.

The existence of a Fourier-Stieltjes algebra augurs well for future applications. In particular, one example suggests an interesting possibility: The algebra of continuous functions on $X$ vanishing at infinity, $C_0(X)$, is the Fourier algebra of a locally compact space $X$. This opens up an entire "dual" approach to the currently exploding subject of non-commutative geometry, which at the moment is regarded more or less exclusively in terms of the associated $C^*$-algebras (not the Fourier-Stieltjes algebras).

As for groups, the Fourier-Stieltjes algebra of a groupoid is the linear span of the positive definite functions and the algebra structure is given by pointwise operations. To provide the Banach space structure, we use $C^*$-algebras attached to $G$, but we use them in a different way from Eymard, and also use $C^*$-algebras associated with the equivalence relation that $G$ induces on $X$.

To describe the various algebras, let us begin with the space $\mathcal{M}_c(G)$ of compactly supported bounded Borel functions on $G$, and its subspace $C_c(G)$. Both are algebras under convolution, which is defined by using the Haar system, and have involutions. If $R$ is the equivalence relation on $X$ induced by $G$, defining $\theta(\gamma) = (r(\gamma), s(\gamma))$ gives a continuous homomorphism of $G$ onto $R$ using the relative product topology on $R$. The quotient topology on $R$ has some advantages: for example, if $\theta$ is one-to-one then $\theta$ is a homeomorphism. ($G$ is said to be *principal*.) Under the quotient topology $R$ is $\sigma$-compact and we can provide it with a Borel measurable Haar system, which allows us to make a convolution $*$-algebra of the space $\mathcal{M}_{\theta c}(R)$ of bounded Borel functions on $R$ that are supported by the image of some compact set in $G$. In Section 1, we show how to make an algebra on $G$ that contains a copy of the space $\mathcal{M}(X)$ of bounded Borel functions on $X$ as well as $\mathcal{M}_c(G)$, and this algebra is denoted by $\mathcal{M}_c(G, X)$. The analog for $R$ is denoted by $\mathcal{M}_{\theta c}(R, X)$. Let $\bar{X}$ denote the one-point compactification of $X$. Then $C(\bar{X}) \subseteq \mathcal{M}(X)$, so $\mathcal{M}_c(G, X)$ contains both $C_c(G)$ and $C(\bar{X})$. The span of these two subalgebras is denoted $C_c(G, \bar{X})$.

If $\omega$ is the universal representation of $G$, then $\omega$ carries each convolution algebra on $G$ to an algebra of operators and thereby provides the convolution algebra with a norm. The closures of the algebras of operators or the completions under the norms are useful in various ways, so we have notation for them: $C^*(G)$ is the completion of $C_c(G)$, $C^*(G, \bar{X})$ is the completion of $C_c(G, \bar{X})$, $M^*(G)$ is the completion of $\mathcal{M}_c(G)$, and $M^*(G, X)$ is the completion of $\mathcal{M}_c(G, X)$. Likewise for $R$ we get $M^*(R)$ and $M^*(R, X)$ from $\mathcal{M}_{\theta c}(R)$ and $\mathcal{M}_{\theta c}(R, X)$. The algebra $\mathcal{B}(G)$ is isomorphic to a Banach algebra of completely bounded operators on $M^*(G)$, but the functions also correspond to completely bounded bimodule mappings from $C^*(G, \bar{X})$ to $M^*(R, X)$ as bimodules over $C(\bar{X})$.

The reader who is unfamiliar with groupoids will find some necessary background material in Section 1. Assuming that background, we define a bounded Borel function $p$ on a locally compact groupoid $G$ with Haar system $\lambda$ to be *positive definite* if

$$\int \int f(\gamma_1) \bar{f}(\gamma_2) p(\gamma_2^{-1} \gamma_1) d\lambda^x(\gamma_1) d\lambda^x(\gamma_2) \geq 0$$

for every $f \in C_c(G)$. The set of these is denoted $\mathcal{P}(G)$ and by definition the set $\mathcal{B}(G)$ is the linear span of $\mathcal{P}(G)$. In both sets two elements that agree except on a negligible set need to be identified, though we find it convenient to indulge in the usual carelessness about maintaining the distinction. The primary result of this paper is

(I) $\mathcal{B}(G)$ is a Banach algebra.

Results needed to prove this are:

(II) Each $p \in \mathcal{P}(G)$ can be represented in terms of a unitary representation of $G$ and a cyclic 'vector' for the representation.

(III) Multiplication by a $b \in \mathcal{B}(G)$ defines a completely bounded operator on $M^*(G)$ whose norm is at least the supremum norm of $b$.

(IV) The set of operators arising from elements of $\mathcal{P}(G)$ is closed in the space of completely bounded operators on $M^*(G)$.

In fact, $\mathcal{B}(G)$ is a Banach algebra of completely bounded operators on $M^*(G)$, and the elements of $\mathcal{P}(G)$ occur as completely positive operators. In order to prove the completeness of $\mathcal{B}(G)$, we introduce an auxiliary groupoid. Let $T_2$ denote the transitive equivalence relation on the two point set $\{1, 2\}$, so that functions on $T_2$ are $2 \times 2$ matrices. Thus functions on $G \times T_2$ can be regarded as $2 \times 2$ matrices of functions on $G$. Then each $b \in \mathcal{B}(G)$ appears as a corner entry of a positive definite function on $G \times T_2$ whose completely bounded norm is the same as that of $b$. Furthermore, such a corner entry is always in $\mathcal{B}(G)$. Combining these facts with the completeness of $\mathcal{P}(G \times T_2)$ is what allows us to finish the proof of completeness of $\mathcal{B}(G)$. This groupoid is used in a similar way in [R4].

The material in this paper can be outlined as follows. Section 1 is devoted to background material on three topics: locally compact groupoids, convolution algebras attached to them, and representations of groupoids and the algebras. Here the reader will find our notation established in the midst of some examples and discussion of results needed later. Section 2 contains some measure theoretic technicalities about Haar systems and choosing of Borel functions in prescribed equivalence classes modulo null sets. Section 3 contains the fundamental results about positive definite functions. We give the definition of 'positive definite function' and establish the connection between such functions and cyclic unitary representations of $G$. In Section 4 we show that multiplication by a positive definite function is a completely positive operator on $M^*(G)$, using the main result of Section 3. Section 4 also includes the proof that a positive definite function gives rise to a completely positive operator from $C^*(G, X)$ to $M^*(R, X)$. All of these operators are bimodule maps over $C(\bar{X})$, the algebra of continuous functions on the one-point compactification of the space of units of $G$. Section 5 contains results about completely bounded bimodule maps. Finally in Section 6 we are able to complete the proof that the linear combinations of positive definite functions constitute a Banach algebra. Section 7 contains some counter-examples.

**Acknowledgements:** We thank Karla Oty for numerous valuable suggestions for improving this paper. We also mention that Jean Renault has obtained similar results in the context of a groupoid with a Haar system and a quasiinvariant measure [R4]. We have used some of his ideas.

## 1. Background on Groupoids

The purpose of this section is to give the reader a source of some essential information about analysis on groupoids needed in this paper. A more thorough background might be useful, and we omit much historical motivation, but we think the bare minimum is here. Some of the information is not in other papers, and we also have to establish notation for this paper.

### Groupoids

Much of our motivation comes from the fact that group actions give rise to groupoids, and that case was important in the development of the subject. However, we want to present a definition that has a different motivation, hoping to make the idea easier to grasp. E. G. Effros suggested this approach.

Start with two sets, $X$ and $G$, and suppose that $X$ is the set of vertices and $G$ the set of edges of a directed graph. If the structure we are about to describe is present, we say that $G$ is a *groupoid on $X$*. Suppose that we have a mapping taking values in $G$ and defined on the set of pairs of edges for which the first edge starts from the vertex where the second edge terminates. (For a groupoid of mappings, we want the operation to be composition and we want the right hand factor to be applied first.) We want the operation to be associative and to have units and inverses.

To describe this in more detail, we use two functions $r$ and $s$ from $G$ onto $X$, such that each $\gamma \in G$ is an edge from $s(\gamma)$ to $r(\gamma)$. Then for $\gamma$ and $\gamma'$ in $G$, the element $\gamma\gamma'$ of $G$ is defined iff $s(\gamma) = r(\gamma')$. We write $G^{(2)} = \{(\gamma, \gamma') \in G \times G : s(\gamma) = r(\gamma')\}$. We also assume there is given a mapping $x \mapsto i_x$ of $X$ into $G$ and an involution $\gamma \mapsto \gamma^{-1}$ on $G$. Then we require the following properties:

(a) (associativity) If $s(\gamma_1) = r(\gamma_2)$, then $s(\gamma_1\gamma_2) = s(\gamma_2)$, and $r(\gamma_1\gamma_2) = r(\gamma_1)$. If, also, $s(\gamma_2) = r(\gamma_3)$, then $(\gamma_1\gamma_2)\gamma_3 = \gamma_1(\gamma_2\gamma_3)$.

(b) (units) If $x \in X$, then $r(i_x) = s(i_x) = x$. If $\gamma \in G$, then $\gamma i_{s(\gamma)} = i_{r(\gamma)}\gamma = \gamma$.

(c) (inverses) $r(\gamma^{-1}) = s(\gamma)$, $s(\gamma^{-1}) = r(\gamma)$, $\gamma\gamma^{-1} = i_{r(\gamma)}$, and $\gamma^{-1}\gamma = i_{s(\gamma)}$

**Remark:** It is usually convenient to identify $x \in X$ with $i_x \in G$, and thus think of $X$ as a subset of $G$. If $\gamma \in G$ and $x = s(\gamma)$, $y = r(\gamma)$, we think of $\gamma$ as mapping $x$ to $y$ and may write $\gamma : x \to y$. (A groupoid is a small category with inverses.)

**Examples:** (a) Suppose a group $H$ acts on a set $X$ (on the left). Set $G = H \times X$, identify $X$ with $\{e\} \times X$, and define $r(h, x) = hx$, $s(h, x) = x$. Then we can define $(h_1, x_1)(h_2, x_2) = (h_1 h_2, x_2)$ if $x_1 = h_2 x_2$, $i_x = (e, x)$ and $(h, x)^{-1} = (h^{-1}, hx)$, to make a groupoid. (Right actions work better for left Haar measures as we see below, and then we have $s(x, h) = xh$, $r(x, h) = x$.)

(b) To make a groupoid from an equivalence relation $R$ on a set $X$, identify $X$ with the diagonal in $X \times X$, define $r(x, y) = x$, $s(x, y) = y$, $(x, y)(y, z) = (x, z)$ and $(x, y)^{-1} = (y, x)$.

(c) Let $X$ be the set of open sets in $\mathbb{R}^n$, and let $G$ be the set of diffeomorphisms between elements of $X$. For $\gamma \in G$, let $s(\gamma)$ be the domain of the mapping and let $r(\gamma)$ be its range. Let the product be function composition and let the inverse be the inverse of functions.

For other examples see [C1, M1, M2, R1, Re1].

Every groupoid determines a natural equivalence relation on its set of units, namely $x \sim y$ iff there is a $\gamma : x \to y$. The equivalence class of $x$ is denoted $[x]$ and is called its *orbit*. As a subset of $X \times X$, this equivalence relation is $R = \{(r(\gamma), s(\gamma)) : \gamma \in G\}$. The function $\theta = (r, s)$ mapping $G$ to $R$ is a groupoid homomorphism and $G$ is called *principal* iff $\theta$ is one–one, i.e., $G$ is isomorphic to an equivalence relation. If $G$ arises from a group action, $G$ is principal iff the action is free (the only element of the group that has any fixed points is the identity).

If $G$ is a groupoid on $X$, and $Y \subseteq X$ is non-empty, we call $r^{-1}(Y) \cap s^{-1}(Y)$ the *restriction of $G$ to $Y$*, and write $G|Y$ for it. In terms of graphs, $G|Y$ is the set of all edges in $G$ that connect points of $Y$. $G|Y$ is a subgroupoid of $G$, and a groupoid on $Y$. For each $x \in X$, $G|\{x\}$ is a group called the *stabilizer* of $x$ or the *isotropy* of $x$.

If $A$ and $B$ are subsets of a groupoid $G$, we define the product $AB$ of the two sets to be $\{\gamma\gamma' : \gamma \in A, \gamma' \in B, r(\gamma') = s(\gamma)\}$. If $A$ has a single element $\gamma_0$, we write $\gamma_0 B$ for $AB$. Thus $YGY = G|Y$ and $xGx = G|\{x\}$ if $Y \subseteq X$ and $x \in X$. We also use the sets $r^{-1}(x) = xG$ and $s^{-1}(x) = Gx$ when $x \in X$.

A groupoid $G$ is a *Borel groupoid* if $G$ has a Borel structure, $X$ is a Borel set when regarded as a subset of $G$, and $r, s, (\ )^{-1}$ and multiplication are Borel functions. We will consider only Borel groupoids which are at least analytic, and then $X = \{\gamma : r(\gamma) = \gamma\}$ is Borel if $r$ is Borel. A groupoid $G$ is *topological* if it has a topology such that $X$ is closed, and $r, s, (\ )^{-1}$ and multiplications are continuous, while $r$ and $s$ are open. Again these properties are not independent. It is necessary for $r$ to be open in order to prove that $AB$ is open whenever $A$ and $B$ are open [R4].

We write $\mathcal{M}(G)$ for the space of bounded Borel measurable functions on $G$, whenever $G$ is a Borel groupoid. If $G$ has a topology in which it is $\sigma$-compact (a countable union of compact sets), we write $\mathcal{M}_c(G)$ for the subspace of $\mathcal{M}(G)$ of functions having compact support.

If $G$ is an analytic Borel groupoid, we say a measure $\mu$ on $G$ is *quasisymmetric* if it has the same null sets as its image $(\mu)^{-1}$ under $(\ )^{-1}$. Thus $\mu$ and $(\mu)^{-1}$ are in the same measure class, and the measure class $[\mu]$ (set of measures with the same null sets as $\mu$) is invariant under $(\ )^{-1}$. For measures on $G$, this global symmetry is just the same as if $G$ were a group.

In the material following this paragraph, we give the definitions for groupoids that extend the notions of invariance and quasiinvariance of measures under translation on a group or under other actions of the group. Because translation on the left by a groupoid element $\gamma$ makes sense only on $s(\gamma)G$, and similarly for right translation, the notions of invariance and quasiinvariance are more complicated for groupoids than for groups.

Following Alain Connes [C1] we say that a *kernel* is a function $\nu$ assigning a $\sigma$-finite (positive) measure $\nu^x$ on $G$ to each $x \in X$, so that these two statements are true:

(a) $\nu^x(G \backslash xG)$ is always 0. One may say that $\nu^x$ is concentrated on $xG$.

(b) If $f \in \mathcal{M}(G)$, and $f \geq 0$, the function $\nu(f) : X \to [0, \infty]$ defined by $\nu(f)(x) = \nu^x(f) = \int f \, d\nu^x$ is Borel.

Given an element $\gamma \in G$, the mapping $\gamma' \mapsto \gamma\gamma'$ is a Borel isomorphism of $s(\gamma)G$ onto $r(\gamma)G$ and thus maps $\nu^{s(\gamma)}$ to a measure $\gamma\nu^{s(\gamma)}$ on $r(\gamma)G$, for every kernel $\nu$. A kernel $\nu$ is called *left invariant* provided $\nu^{r(\gamma)} = \gamma\nu^{s(\gamma)}$ for all $\gamma \in G$. It is called *(left) quasiinvariant* if $\nu^{r(\gamma)}$ and $\gamma\nu^{s(\gamma)}$ are equivalent for all $\gamma \in G$.

A left invariant kernel, $\lambda$, on a Borel groupoid $G$ is called a *Borel Haar system*. Then defining $\lambda_x$ to be the image of $\lambda^x$ under inversion produces a right Borel Haar system. A Borel Haar system $\lambda$ on a locally compact groupoid is called a *Haar system* if $\text{supp}(\lambda^x)$ is always $xG$ and $\lambda(f) \in C_c(X)$ for each $f \in C_c(G)$. In particular, each $\lambda^x$ is a Radon measure. For discussions of Haar systems, see [H2, Re1, S1, S2].

When $\lambda$ is a Haar system, it can be convenient to have a left quasiinvariant kernel $\lambda_1$ consisting of probability measures equivalent to the measures $\lambda^x$. It is not difficult to show that there is a continuous, strictly positive, function $f$ on $G$ such that for every $x \in X$, $\int f \, d\lambda^x = 1$. We choose one such $f$ and write $\lambda_1^x$ for the measure $f\lambda^x$. We also write $\mu_1^x$ for the probability measure $s(\lambda_1^x)$ on $X$; these measures also depend on $x$ in a Borel way.

If $\lambda$ is a Borel Haar system on a Borel groupoid $G$ and $\mu$ is a probability measure on $X$, we can form a measure $\nu = \int \lambda^x d\mu(x)$: $\int f d\nu = \int \int f(\gamma) d\lambda^x(\gamma) d\mu(x)$. We often write $\lambda^\mu$ for this measure $\nu$. Suppose that $G = X \times H$, where $X$ is a right $H$-space, and give $G$ the groupoid structure that comes from the group action. Let $\lambda$ be a left Haar measure on $H$. For each $x \in X$, let $\epsilon^x$ be the point mass at $x$, and define $\lambda^x = \epsilon^x \times \lambda$, to get a Borel left Haar system. If $\mu$ is a $\sigma$-finite measure on $X$ for this groupoid, then $\nu = \lambda^\mu = \mu \times \lambda$ and the class $[\nu]$ is symmetric iff $\mu$ is quasiinvariant under the group action, i.e., for every Borel set $E \subseteq X$ and every group element $h$, $\mu(E) = 0$ iff $\mu(Eh) = 0$. This follows from [R1, Theorem 4.3] and the fact that if $\mu$ is quasiinvariant under almost all elements of the group, then it is quasiinvariant. Hence, on a general Borel groupoid with Borel Haar system $\lambda$, a $\sigma$-finite measure $\mu$ on $X$ is called *quasiinvariant* iff $\lambda^\mu$ is quasisymmetric. In that case, a result of Peter Hahn [H1, Corollary 3.14], combined with [R3, Theorem 3.2], shows that there is a Borel homomorphism, $\Delta_\mu$, of $G$ to the multiplicative positive real numbers such that

$$\Delta_\mu = \frac{d\lambda^\mu}{d(\lambda^\mu)^{-1}}.$$

This homomorphism is called the *modular function* by analogy with locally compact groups. If $\mu$ is quasiinvariant, and $Y$ is a $\mu$-conull Borel set in $X$, the restriction $G|Y$ is called *inessential*.

We often refer to the set of all quasiinvariant $\sigma$-finite measures on $X$, and will denote that set by $\mathcal{Q}$. We say a Borel set $N \subseteq X$ is $\mathcal{Q}$-*null* provided $\mu(N) = 0$ for every $\mu \in \mathcal{Q}$. It follows from the existence, and uniqueness up to equivalence, of a quasiinvariant $\sigma$-finite measure on each orbit [Re1, Prop. 1.3.6 and Prop. 1.3.8] that $N$ is $\mathcal{Q}$-null iff $\lambda^x(GN)$ is always 0. The measures $\mu_1^x$ introduced above are in this class, and any measure in $\mathcal{Q}$ equivalent to such a measure is called *transitive* because it is concentrated on a single orbit. For a Borel set $N \subseteq G$, we say $N$ is $\lambda^\mathcal{Q}$-null iff $\lambda^\mu(N) = 0$ whenever $\mu \in \mathcal{Q}$. A function $f$ on $X$ is $\mathcal{Q}$-*essentially bounded* iff the restriction of $f$ to the complement of some $\mathcal{Q}$-null set is bounded, and then $\|f\|_\infty$ is defined to be the smallest element of $\{B : |f| \leq B \ \mu$-amost everywhere for every $\mu \in \mathcal{Q}\}$. The space of $\mathcal{Q}$-essentially bounded functions on $X$ will be denoted by $L^\infty(\mathcal{Q})$. A similar definition is used for the space $L^\infty(\lambda^\mathcal{Q})$ of $\lambda^\mathcal{Q}$-essentially bounded functions on $G$, except that the measures $\lambda^\mu$ are used.

**Examples**: (a) If $G = X \times H$, where $X$ and $H$ are locally compact and $H$ is a group, let $\epsilon^x$ denote the unit point mass at $x$ for $x \in X$ and let $\lambda$ be a left Haar measure on $H$. Then $\lambda^x = \epsilon^x \times \lambda$ defines a Haar system for $G$.

(b) If $E$ is an analytic equivalence relation on $X$ and each equivalence class is countable, we can let $\lambda^x$ be counting measure on $\{x\} \times [x]$ to get a left invariant system of measures.

(c) Here is an example of a locally compact groupoid that has a Borel Haar system but no Haar system. Let $G = [0, 1/2] \times \{0\} \cup [1/2, 1] \times \mathbb{Z}/2$. This is a field of groups. To get a Borel Haar system, we can make each $\lambda^x$ a multiple of the Haar measure on $\{0\}$ or $\mathbb{Z}/2$. Then $\lambda^{1/2}(\{(1/2,0)\}) = \lambda^{1/2}(\{(1/2,1)\}) > 0$, and if we let $f$ be the characteristic function of $[1/2, 1] \times \{1\}$, then the function $\lambda(f)$ has a jump at $1/2$. We could easily change to another locally compact topology on this $G$ and get a Haar system. In general, it may be necessary to change the topology on $G$ and pass to an inessential restriction in order to get a Haar system.

**Convolution**

We use several convolution algebras in this paper, and will introduce them here. There are two basic convolutions, a convolution of functions that can be defined in the presence of a Borel Haar system, and a convolution of kernels that does not depend on any such system. If the groupoid is locally compact and the Haar system is continuous, then $C_c(G)$ is an algebra under the convolution of functions. We will see that convolution of functions can be subsumed under convolution of kernels by replacing each function by the kernel obtained by multiplying the Haar system by the function.

First, let $G$ be a Borel groupoid with a Borel Haar system $\lambda$. If $f, g$ are non-negative Borel functions on $G$, then $\int f(\gamma_1)g(\gamma_2)d\lambda^{r(\gamma_1)}(\gamma_2)$ is a Borel function of $\gamma_1$, so by taking linear combinations and monotone limits we see that whenever $F$ is a non-negative Borel function on $G \times G$ the integral $\int F(\gamma_1, \gamma_2) d\lambda^{r(\gamma_1)}(\gamma_2)$ depends on $\gamma_1$ in a Borel manner. Then for non-negative $f, g \in \mathcal{M}(G)$, we can let $F(\gamma_1, \gamma_2) = f(\gamma_1)g(\gamma_1^{-1}\gamma_2)$ when $r(\gamma_2) = r(\gamma_1)$ and $F(\gamma_1, \gamma_2) = 0$ otherwise, and see that $\int f(\gamma_1)g(\gamma_1^{-1}\gamma_2) d\lambda^{r(\gamma_2)}(\gamma_1)$ is a Borel function of $\gamma_2$. Denote this function by $f * g$, provided that it is always finite valued. Then $f * g \in \mathcal{M}(G)$. The function $f * g$ is called the *convolution* of $f$ and $g$. Convolution can be extended to more general functions using linearity.

Define the space $I_r(G, \lambda)$ to be $\{f \in \mathcal{M}(G) : \lambda(|f|) \text{ is bounded}\}$, and give it a norm by letting $\|f\|_{I,r}$ be the sup norm of the Borel function $\lambda(|f|)$. As proved on page 51 of [Re1], $I_r(G, \lambda)$ is closed under convolution and the norm $\|\cdot\|_{I,r}$ is an algebra norm. We can define an involution on $\mathcal{M}(G)$ as in [Re1] by letting $f^\flat(\gamma) = \bar{f}(\gamma^{-1})$ for $f \in \mathcal{M}(G)$, $\gamma \in G$. If we set $I(G, \lambda) = I_r(G, \lambda) \cap (I_r(G, \lambda))^\flat$, then we can define $\|f\|_I$ to be the maximum of $\|f\|_{I,r}$ and $\|f^\flat\|_{I,r}$ for $f \in I(G, \lambda)$, obtaining a normed algebra on which the involution is an isometry.

If $G$ is locally compact and $\lambda$ is a Haar system, then $C_c(G)$ is a $*$-subalgebra of $I(G, \lambda)$. In the inductive limit topology, $C_c(G)$ is a topological algebra [Re1, page 48].

The second kind of convolution can be introduced after the objects are defined: A *complex kernel* is a function $\nu$ assigning a complex measure $\nu^x$ on $G$ so that

(a) $\nu^x$ is always concentrated on $xG$.

(b) if $f \in \mathcal{M}(G)$, the function $\nu(f)$ taking $x \in X$ to $\nu^x(f)$ is Borel.

We define $K(G)$ to be the space of *bounded* complex kernels on $G$, i.e., those for which the total variation of $\nu^x$ is a bounded function of $x$.

If $\gamma \in G$ and $\nu \in K(G)$ we can map $\nu^{s(\gamma)}$ to a measure on $r(\gamma)G$, via left translation by $\gamma$, as we do in defining Haar systems. Denote this measure by $\gamma\nu^{s(\gamma)}$. If $\nu_1, \nu_2 \in K(G)$ we can define the convolution $\nu = \nu_1 * \nu_2$ by $\nu^x = \int \gamma\nu_2^{s(\gamma)} d\nu_1^x(\gamma)$. We can also define an action of $K(G)$ on $I_r(G, \lambda)$ as follows. If $\nu \in K(G)$, $f \in I_r(G, \lambda)$ and $\gamma' \in G$ set

$$L(\nu)f(\gamma') = \int f(\gamma^{-1}\gamma') d\nu^{r(\gamma')}(\gamma).$$

It is not difficult to verify that $L(\nu)$ is a bounded operator whose norm is at most the essential supremum of the total variation norms of the signed measures $\nu^x$. If $\nu_1$ and $\nu_2$ are in $K(G)$ and $f \in I_r(G, \lambda)$, then we can calculate as follows:

$$\begin{aligned}(L(\nu_1)(L(\nu_2)f))(\gamma) &= \int (L(\nu_2)f)(\gamma_1^{-1}\gamma) d\nu_1^{r(\gamma)}(\gamma_1) \\ &= \int\int f(\gamma_2^{-1}\gamma_1^{-1}\gamma) d\nu_2^{s(\gamma_1)}(\gamma_2) d\nu_1^{r(\gamma)}(\gamma_1) \\ &= \int\int f(\gamma_2^{-1}\gamma) d\nu_2^{s(\gamma_1)}(\gamma_1^{-1}\gamma_2) d\nu_1^{r(\gamma)}(\gamma_1) \\ &= \int f(\gamma_2^{-1}\gamma) d(\nu_1 * \nu_2)^{r(\gamma)}(\gamma_2),\end{aligned}$$

showing that $L$ takes convolution to composition of operators. Since $L$ is faithful, $K(G)$ is an algebra under convolution. If $f, g \in I_r(G, \lambda)$ it is not difficult to verify that $f\lambda \in K(G)$ and $L(f\lambda)g = f * g$:

$$L(f\lambda)g(\gamma) = \int g(\gamma_1^{-1}\gamma)f(\gamma_1) d\lambda^{r(\gamma)}(\gamma_1).$$

Since $L$ is faithful and convolution is associative, it follows that $f\lambda * g\lambda = (f * g)\lambda$. Thus $I_r(G, \lambda)\lambda = \{f\lambda : f \in I_r(G, \lambda)\}$ is a subalgebra of $K(G)$ isomorphic to $I_r(G, \lambda)$. If $G$ is locally compact and has a Haar system $\lambda$, the calculations just made also show that $C_c(G)\lambda$ is a subalgebra of $K(G)$ isomorphic to $C_c(G)$.

Next we want to enlarge $C_c(G)\lambda$ to a subalgebra of $K(G)$ that contains a copy of $C_c(X)$. We denote the one-point compactification of $X$ by $\bar{X}$. The mapping $f \to f|X$ takes $C(\bar{X})$ one-one onto the algebra of continuous functions on $X$ that have a limit at infinity. We identify $C(\bar{X})$ with that subalgebra of $C(X)$ but continue to write $C(\bar{X})$. Notice that there is also a subalgebra of $K(G)$ isomorphic to $C(\bar{X})$, obtained as follows. First define $\epsilon$ to be the kernel that assigns the point mass at $x$ to each $x \in X$, which we denote by $\epsilon^x$ as above. Next notice that $K(G)$ is closed under multiplication by any bounded Borel function on $G$, so if $h \in \mathcal{M}(X)$ and $\nu \in K(G)$, we can define $h\nu$ to be $(h \circ r)\nu$, and $\nu h = (h \circ s)\nu$. (These agree with the naturally defined left and right multiplication of $\mathcal{M}(X)$ on $I_r(G, \lambda)$ when the latter is regarded as a space of kernels.) Then $\mathcal{M}(X)\epsilon$ is a subalgebra of $K(G)$ isomorphic to $\mathcal{M}(X)$, and that algebra includes $C(\bar{X})\epsilon$, which is isomorphic to $C(\bar{X})$.

If we write $C_c(G, \bar{X})$ for the sum of $C(\bar{X})\epsilon$ and $C_c(G)\lambda$ as subspaces of $K(G)$, it can be seen that $C_c(G, \bar{X})$ is a subalgebra. Also the involution on $C_c(G)$ extends in a natural way to $C_c(G, \bar{X})$. We need the algebra $C_c(G, \bar{X})$ because it generates a $C^*$-algebra that contains $C(\bar{X})$ as a subalgebra, enabling us to apply a result on completely bounded bimodule mappings in Section 5.

On the other hand, the algebra $C_c(G)$ has an approximate unit [M-R-W]. In order to state the existence theorem, we need to introduce some of their terminology. They call a set $L$ in $G$ *r-relatively compact* if $KL$ is relatively compact for every compact set $K \subseteq X$. There exists a decreasing sequence $U_1, U_2, \ldots$ of open r-relatively compact sets whose intersection is $X$. There also exists an increasing sequence of compact sets in $X$, $K_1, K_2, \ldots$ whose interiors exhaust $X$. These come from the second countability of $G$, and they allow us to make a sequence that is an approximate unit (instead of a more general net). We call a function $f$ in $C_c(G)$ *symmetric* if $f^\flat = f$.

**Theorem** 1.1. [M-R-W] There is a sequence $e_1, e_2, \ldots$ of symmetric functions in $C_c^+(G)$ such that for each $n$ we have

(i) $\mathrm{supp}(e_n) \subseteq U_n$, and

(ii) $\int e(\gamma) \, d\lambda^x(\gamma) \geq 1 - n^{-1}$ for $x \in K_n$, and $\leq 1$ for all $x \in X$.

Such a sequence is a two-sided approximate unit for $C_c(G)$ in its inductive limit topology, i.e., for uniform convergence on compact sets.

**Representations**

A (unitary) *representation* of a locally compact groupoid $G$ is given by a Hilbert $G$-bundle $\mathcal{K}$ over $X$, the unit space of $G$; this means we have two functions that have some properties:

(a) a Hilbert space $\mathcal{K}(x)$ for each $x$. We form $\Gamma_\mathcal{K} = \{(x, v) : x \in X, v \in \mathcal{K}(x)\}$, called the *graph of $\mathcal{K}$*, and require that $\Gamma_\mathcal{K}$ have a standard Borel structure such that the projection onto $X$ is Borel and there is a countable set of Borel sections of $\Gamma_\mathcal{K}$ such that for each $x$ the set of their values at $x$ is dense in $\mathcal{K}(x)$.

(b) a Borel homomorphism $\pi$ of $G$ into the unitary groupoid of the bundle $\mathcal{K}$, i.e., for each $\gamma$, $\pi(\gamma) : \mathcal{K}(s(\gamma)) \to \mathcal{K}(r(\gamma))$ is unitary, and $\pi$ is a Borel function [Re1, page 52].

This can also be said as follows: $(\mathcal{K}, \pi)$ is a Borel functor on $G$ taking values in the category of Hilbert spaces.

Given a representation $\pi$ of $G$, and a measure $\mu \in \mathcal{Q}$, we can obtain from them a $*$-representation of $\mathcal{M}_c(G)$. Before describing the representation, we need another item of notation. We will write $\nu = \lambda^\mu = \int \lambda^x d\mu(x)$. Then we take $\Delta = \Delta_\mu$ as above and define $\nu_0 = \Delta^{-\frac{1}{2}}\nu$, obtaining a symmetric measure. Next we make a Hilbert space, $L^2(\mu; \mathcal{K})$, of square integrable sections of $\mathcal{K}$. For $f \in \mathcal{M}_c(G)$ we define $\pi^\mu(f)$ on $L^2(\mu; \mathcal{K})$ by setting

$$(\pi^\mu(f)\xi \mid \eta) = \int f(\gamma)(\pi(\gamma)\xi \circ s(\gamma) \mid \eta \circ r(\gamma)) \, d\nu_0(\gamma)$$

for $\xi, \eta \in L^2(\mu; \mathcal{K})$. Then $\pi^\mu$ is a $*$-representation of $\mathcal{M}_c(G)$ with $\|\pi^\mu(f)\| \leq \|f\|_{I,\mu}$ [Re1, pages 52,53], so its restriction to $C_c(G)$ has the same property. We denote the restriction by the same symbol, depending on context to distinguish the two. Later we will also use another method of integrating a unitary representation of $G$, one that is due to Hahn [H2] and does not use the symmetrized measure.

It can be convenient to choose $\mu$ to be finite, say a probability measure, so we need to know that $\mu' \sim \mu$

implies $\pi^{\mu'}$ is unitarily equivalent to $\pi^\mu$. To prove this implication, take $\rho$ to be a positive Borel function whose square is the Radon-Nikodym derivative of $\mu'$ with respect to $\mu$. Then

$$\rho^2 \circ r = \frac{d\lambda^{\mu'}}{d\lambda^\mu}$$

and

$$\rho^2 \circ s = \frac{d(\lambda^{\mu'})^{-1}}{d(\lambda^\mu)^{-1}}$$

so

$$(\rho^2 \circ r)\Delta_\mu = (\rho^2 \circ s)\Delta_{\mu'}.$$

Hence we can define $V : L^2(\mu', \mathcal{K}) \to L^2(\mu, \mathcal{K})$ by $V\xi = \rho\xi$ to get the necessary unitary equivalence. To see that it is indeed an intertwining operator, compute to see that the inner products are equal: $(\pi^\mu(f)V\xi \mid V\eta) = (\pi^{\mu'}(f)\xi \mid \eta)$.

It is natural to ask whether every continuous representation of $C_c(G)$ can be obtained by integrating a unitary representation of $G$, as is true for groups. An affirmative answer to this question was provided by an ingenious argument due to Jean Renault [Re2], and it follows that every representation of $\mathcal{M}_c(G)$ bounded by $\|\ \|_I$ can be obtained by integrating a unitary representation of $G$. Another discussion of this result is in [Mu]. Renault's theorem is:

**Theorem 1.2**. Let $G$ be a locally compact groupoid that has a Haar system, and let $\mathcal{H}_0$ be a dense subspace of a (separable) Hilbert space $\mathcal{H}$. Suppose that $L$ is a representation of $C_c(G)$ by operators on $\mathcal{H}_0$ such that

(a)  $L$ is non-degenerate;

(b)  $L$ is continuous in the sense that for every pair of vectors $\xi, \eta \in \mathcal{H}_0$, the linear functional $L_{\xi,\eta}$ defined by $L_{\xi,\eta}(f) = (L(f)\xi \mid \eta)$ is continuous relative to the inductive limit topology on $C_c(G)$;

(c)  $L$ preserves the involution, i.e., $(\xi \mid L(f^\flat)\eta) = (L(f)\xi \mid \eta)$ for $\xi, \eta \in \mathcal{H}_0$ and $f \in C_c(G)$.

Then the operators $L(f)$ are bounded. The representation of $C_c(G)$ on $\mathcal{H}$ obtained from $L$ is equivalent to one obtained by integrating a unitary representation of $G$ using a probability measure $\mu \in \mathcal{Q}$. In particular, $L$ is continuous relative to $\|\ \|_I$.

In [Re1], Renault defined a norm on $C_c(G)$ by $\|f\| = \sup\{\|L(f)\| : L$ is a bounded representation of $C_c(G)\}$. Theorem 1.2 shows that we could get the same norm by using the representations $\pi^\mu$ in place of the $L$'s. The completion of $C_c(G)$ with respect to the norm just defined is a $C^*$-algebra denoted $C^*(G)$. Every positive linear functional of norm one on a $C^*$-algebra gives rise to a representation of the algebra and a cyclic vector in the Hilbert space of the representation. The direct sum of all these cyclic representations is called the *universal representation* of the $C^*$-algebra. We will denote this representation by $\omega$. For $C^*(G)$, we know that every one of the cyclic representations is of the form $\pi^\mu$, so $\omega$ can also be regarded as a representation of $\mathcal{M}_c(G)$. We will write $M^*(G)$ for the operator norm closure of $\omega(\mathcal{M}_c(G))$. Since $\omega$ is an isomorphism on $C^*(G)$, we can regard $C^*(G)$ as a subalgebra of $M^*(G)$. We will also refer to $\omega$ as the *universal representation* of $G$ itself.

In proving that $L$ can be obtained by integration, Renault shows that there is a representation of $C_c(X)$, say $\phi$, associated with $L$ such that for $f \in C_c(G)$ and $h \in C_c(X)$ we have

$$L((h \circ r)f) = \phi(h)L(f)$$

and

$$L(f(h \circ s)) = L(f)\phi(h).$$

Then $\phi$ extends in the obvious way to a unital representation of $C(\bar{X})$ and can be used to extend $L$ to a representation of $C_c(G, \bar{X})$:

$$L(f\lambda + g\epsilon) = L(f) + \phi(g).$$

The reader can verify, easily, that this defines a unital representation of $C_c(G, \bar{X})$. We extend $\omega$ to $C(G, \bar{X})$ in this way, and also to $\mathcal{M}_c(G, \bar{X})$. Then we define $C^*(G, \bar{X})$ to be the operator norm closure of $\omega(C_c(G, \bar{X}))$ and $M^*(G, X)$ to be the closure of $\omega(\mathcal{M}_c(G, X))$.

For some computations we need another norm [H2]. Let $\mu \in \mathcal{Q}$, let $f \in \mathcal{M}_c(G)$ and define

$$\|f\|_{II,\mu} = \sup\{\int |f(\gamma)g \circ r(\gamma)h \circ s(\gamma)|\Delta_\mu(\gamma)^{-1/2}d\lambda^\mu(\gamma)\},$$

the supremum being taken over unit vectors $g, h \in L^2(\mu)$. Then define $\|f\|_{II}$ to be $\sup\{\|f\|_{II,\mu} : \mu \in \mathcal{Q}\}$. Three facts about this norm should be mentioned. The first is that if $\pi$ is a unitary representation of $G$, then $\|\pi^\mu(f)\| \leq \|f\|_{II,\mu}$. Thus $\|\omega(f)\| \leq \|f\|_{II}$, because $\|\omega(f)\| = \sup\{\|\pi^\mu(f)\| : \pi$ is a unitary representation and $\mu \in \mathcal{Q}\}$. Next, if $\pi$ is the one dimensional trivial representation and $f \geq 0$ then $\|\pi^\mu(f)\| = \|f\|_{II,\mu}$. It follows that if $0 \leq f \in \mathcal{M}_c(G)$ then

$$\|\omega(f)\| = \|f\|_{II}.$$

A third fact is this: if $b \in L^\infty(\lambda^\mathcal{Q})$ and $f \in \mathcal{M}_c(G)$, then for any $\mu \in \mathcal{Q}$ we have

$$\|bf\|_{II,\mu} \leq \|b\|_\infty \|f\|_{II,\mu},$$

so

$$\|bf\|_{II} \leq \|b\|_\infty \|f\|_{II}.$$

**Lemma 1.3.** If $0 \leq f \in \mathcal{M}_c(G)$ and $b \in L^\infty(\lambda^\mathcal{Q})$, then $\|\omega(bf)\| \leq \|b\|_\infty \|\omega(f)\|$.

**Proof.** Using the three properties of $\|\ \|_{II,\mu}$ mentioned just above, we have

$$\begin{aligned}\|\omega(bf)\| &\leq \sup\{\|bf\|_{II,\mu} : \mu \in \mathcal{Q}\} \\ &\leq \sup\{\|b\|_\infty \|f\|_{II,\mu} : \mu \in \mathcal{Q}\} \\ &= \|b\|_\infty \|\omega(f)\|.\end{aligned}$$

## 2. Measure Theoretic Preparation

A basic lemma is needed for our construction of positive definite functions from completely positive maps in Section 4. After proving that lemma, we also need to prepare some detailed information about Haar systems on locally compact groupoids and how they relate to Borel Haar systems on the associated equivalence relations. Most of that information comes from [Re3].

As preparation for the proof of the lemma in question, we recall a basic fact about measures and function spaces. Suppose that $(X, \mathcal{B})$ is a set with a $\sigma$-algebra and that $\mathcal{A}$ is a subalgebra of $\mathcal{B}$ that generates $\mathcal{B}$ as a $\sigma$-algebra. Let $\mu$ be any finite measure defined on $\mathcal{B}$. The measure of the symmetric difference between two sets is the same as the distance between their characteristic functions in $L^1(\mu)$, and hence provides a (pseudo) metric on $\mathcal{B}$. The closure of $\mathcal{A}$ in $\mathcal{B}$ is a $\sigma$-algebra that contains $\mathcal{A}$ and hence is $\mathcal{B}$. For us, it is important that the fact of density is independent of $\mu$. This implies similar properties for the set $S(\mathcal{A})$, our notation for the set of linear combinations of characteristic functions of sets in $\mathcal{A}$ using coefficients from $\mathbb{Q}[i]$, which is $\mathbb{Q}$ with $\sqrt{-1}$ adjoined. By looking first at simple functions, it is easy to show that $S(\mathcal{A})$ is always dense in $L^1(\mu)$. In the same way, we see that for any $f \in L^1(\mu)$,

$$\|f\|_1 = \sup\{|\int f\varphi\, d\mu| : \varphi \in S(\mathcal{A}) \text{ and } |\varphi| \leq 1\},$$

which is a supremum indexed by a family independent of $\mu$. When $\mathcal{A}$ can be taken to be countable, as is the case when $X$ is a standard Borel space, these facts are particularly useful.

A similar situation arises if $X$ is locally compact. In that case, there is a countable dense subset $\mathcal{S}(X)$ of $C_c(X)$ that is an algebra over $\mathbb{Q}[i]$, and any such $\mathcal{S}(X)$ is dense in $L^1(\mu)$ for every finite measure $\mu$ on $X$.

The next lemma is a generalization of the fact that for two measure spaces, functions on the product and functions from one measure space to the functions on the other can be identified. In our setting, the measure on the image space must be allowed to vary.

**Lemma 2.1.** Let $X$ and $Y$ be standard Borel spaces and let $x \mapsto \nu^x$ be a Borel function from $X$ to finite Borel measures on $Y$. Suppose that $f$ is a function on $X$ selecting an element $f(x)$ of $L^1(\nu^x)$ for each $x \in X$ so that the function $x \mapsto f(x)\nu^x$ is Borel, taking values in the space of complex valued Borel measures. Then there is a Borel function $F$ on $X \times Y$ such that for every $x \in X$ the function $F(x, \cdot)$ is integrable relative to $\nu^x$ and in the class $f(x)$. The function $F$ can be chosen so that if $f(x) \in L^\infty(\nu^x)$ then $F(x, \cdot)$ is bounded by $\|f(x)\|_\infty$. It is possible to choose $F$ meeting those conditions and so that if $\nu^x = \nu^{x'}$ and $f(x) = f(x')$ then $F(x, \cdot) = F(x', \cdot)$ (everywhere on $Y$).

**Remark.** If $Y$ is taken to be locally compact, the proof given below can be modified by replacing $S(\mathcal{A})$ by $\mathcal{S}(Y)$.

**Proof.** For the proof we must have a way, that does not depend on $x$ directly, to choose representatives of classes approximating $f(x)$. For this we choose first a countable algebra, $\mathcal{A}$, of Borel sets in $Y$ that generates the $\sigma$-algebra of Borel sets, so we can use the facts mentioned before the statement of the lemma. List $S(\mathcal{A})$ as a sequence, $s_1, s_2, \ldots$. For convenience, let us write $x \sim x'$ to mean that $\nu^x = \nu^{x'}$ and $f(x) = f(x')$, and say that such points are *equivalent*.

Now we are ready to describe the basic step which will be used repeatedly in the proof. If $\epsilon > 0$ and $x \in X$ define $j(x, \epsilon)$ to be the least element of $\{i : \|f(x) - s_i\|_{L^1(\nu^x)} < \epsilon\}$. It is clear that $j(\cdot, \epsilon)$ takes the same value at equivalent points of $X$, and we will show that $j(\cdot, \epsilon)$ is a Borel function. This will follow if we can show that for each bounded Borel function $h$ on $Y$, $\{x : \|f(x) - h\|_{L^1(\nu^x)} < \epsilon\}$ is a Borel set. We can get that from the fact that norms can be computed as suprema, because for each $\varphi \in S(\mathcal{A})$, $\int (f(x) - h)\varphi \, d\nu^x$ is a Borel function of $x$ and hence so is its absolute value.

If we define $g(x) = s_{j(x,\epsilon)}$ (as an element of $L^1(\nu^x)$) and $G(x, y) = s_{j(x,\epsilon)}(y)$, then $g(x) = g(x')$ and $G(x, \cdot) = G(x', \cdot)$ (everywhere on $Y$) whenever $x \sim x'$. Also, both these functions are Borel.

Apply this process first to $f$ with $\epsilon = 2^{-1}$ to obtain $G_1$ and $g_1$. Then apply it to $f - g_1$ with $\epsilon = 2^{-2}$ to obtain $G_2$ and $g_2$, etc. For each $n$ the value of the function $f - (g_1 + \cdots + g_n)$ at a point $x$ is an element of $L^1(\nu^x)$ having norm $< 2^{-n}$. Thus for $n \geq 2$, $\|g_n(x)\|_1 < 3(2^{-n})$. It follows that for each $x$ the sum $\sum_{n \geq 1} |G_n(x, y)|$ is finite for almost all $y$. Inductively, we see that $G_n(x, \cdot) = G_n(x', \cdot)$ if $x \sim x'$. The set $N = \{(x, y) \in X \times Y : \sum_{n \geq 1} |G_n(x, y)| = \infty\}$ is a Borel set in $X \times Y$ and the slices of $N$ over $x$ and $x'$ are the same set if $x \sim x'$. Now change each $G_n$ to be 0 on $N$. Then the sum is always finite and we still have $G_n(x, \cdot) = G_n(x', \cdot)$ if $x \sim x'$.

Define $F(x, y) = \sum_{n \geq 1} G_n(x, y)$. Then $F$ is Borel and satisfies the first and last conclusions of the theorem. Thus the slice of the Borel set $\{(x, y) : |F(x, y)| > \|f(x)\|_\infty\}$ over every point of $X$ is of measure 0 and the slices of this set are the same over equivalent points of $X$. Change $F$ to be 0 on that set, and all the desired conditions are satisfied.

Now we are going to present some results on the fine structure of the Haar system, as developed by J. Renault in Section 1 of [Re3]. Renault decomposes the Haar system $\lambda$ over a Borel Haar system $\alpha$ on $R$, by studying the action of $G$ on a special group bundle, and we summarize the results here. Recall that the isotropy group bundle of $G$, denoted by $G'$, is defined to be $\{\gamma \in G : r(\gamma) = s(\gamma)\} = \bigcup \{xGx : x \in X\}$. This is closed in $G$ and hence locally compact, so the space of closed subsets of $G'$ is a compact space in the Fell topology [Fell]. Let $\Sigma^{(0)}$ be the space of closed subgroups of the fibers in $G'$, which is a closed subset of the space of closed subsets. Then the set $\Sigma = \{(H, \gamma) \in \Sigma^{(0)} \times G' \mid \gamma \in H\}$ is called the *canonical group bundle* of $\Sigma^{(0)}$. $G$ acts on $\Sigma$ and on $\Sigma^{(0)}$ by conjugation: if $(H_1, \gamma_1) \in \Sigma$, $\gamma \in G$, and $s(\gamma_1) = r(\gamma)$, then

$$(H_1, \gamma_1)\gamma = (\gamma^{-1} H_1 \gamma, \gamma^{-1} \gamma_1 \gamma),$$

while if $H \in \Sigma^{(0)}$, say $H \subseteq xGx$, and $r(\gamma) = x$, then $H \cdot \gamma = \gamma^{-1} H \gamma$. We want to make a Borel choice of Haar measures on the groups $xGx$. One way to do this is to choose a continuous function $F_0$ on $G$ that is non-negative, 1 at each $x \in X$ and has compact support on each $xG$. Then for each $x \in X$ choose a left Haar measure $\beta^x$ on $xGx$ so the integral of $F_0$ with respect to $\beta^x$ is 1. Likewise, choose a function $F$ on $\Sigma$ that is non-negative, 1 at each point $(H, e)$, and has support that intersects every $\{H\} \times H$ in a compact set, and make a similar choice of Haar system on $\Sigma$, $\beta^H$.

Form the groupoid $\Sigma^{(0)} * G = \{(H, \gamma) : s(H) = r(\gamma)\}$ arising from the action of $G$ on $\Sigma^{(0)}$. Then the essential uniqueness of Haar measures guarantees the existence of a 1-cocycle, $\delta$, on $\Sigma^{(0)} * G$ so that for every $(H, \gamma) \in \Sigma^{(0)} * G$ we have
$$\gamma^{-1}\beta^H \gamma = \delta(H, \gamma)^{-1} \beta^{\gamma^{-1}H\gamma}.$$
Renault proves that $\delta$ is continuous. The cohomology class of $\delta$ is determined by $G$, and Renault calls it the *isotropy modulus function of $G$*.

To shorten some formulas in this context, we write $G(x)$ for $xGx$. Renault defines $\delta(\gamma) = \delta(G(r(\gamma)), \gamma)$ to get a 1-cocycle, also called $\delta$, on $G$ such that for every $x \in X$, $\delta|xGx$ is the modular function for $\beta^x$. The pre-image in $\Sigma^{(0)} * G$ of a compact set in $G$ is compact, so $\delta$ and $\delta^{-1} = 1/\delta$ are bounded on compact sets in $G$. Renault defines $\beta^x_y = \gamma \beta^y$ if $\gamma \in xGy$. If $\gamma'$ is another element of $xGy$, then $\gamma^{-1}\gamma' \in yGy$, and since $\beta^y$ is a left Haar measure on $yGy$, it follows that $\beta^x_y$ is independent of the choice of $\gamma$. With this apparatus in place, it is possible to describe a decomposition of the Haar system $\lambda$ for $G$ over the equivalence relation $R = \{(r(\gamma), s(\gamma)) : \gamma \in G\}$. This $R$ is the image of $G$ under the homomorphism $\theta (= (r,s))$, so it is a $\sigma$-compact groupoid. Furthermore, there is a unique Borel Haar system $\alpha$ for $R$ with the property that for every $x \in X$ we have
$$\lambda^x = \int \beta^z_y \, d\alpha^x(z, y).$$
Now suppose that $\mu \in \mathcal{Q}$ so that we can form $\alpha^\mu$ and $\lambda^\mu$, getting quasisymmetric measures. If $\underline{\Delta} = d\alpha^\mu / d(\alpha^\mu)^{-1}$ then $\delta \underline{\Delta} \circ \theta$ will serve as $d\lambda^\mu / d(\lambda^\mu)^{-1}$, i.e., we can always take $\Delta_\mu = \delta \underline{\Delta}_\mu \circ \theta$. We shall see that sometimes $\underline{\Delta}_\mu = 1$ so $\Delta_\mu = \delta$.

For each $x$, the measure $\alpha^x$ is concentrated on $\{x\} \times [x]$ so there is a measure $\mu^x$ on $[x]$ such that $\alpha^x = \epsilon^x \times \mu^x$, where $\epsilon^x$ is the unit point mass at $x \in X \subseteq G$. Since $\alpha$ is a Haar system, we have $\mu^x = \mu^y$ if $x \sim y$, and the function $x \mapsto \mu^x$ is Borel. If we take $\mu'$ to be the measure $\mu^z$ for some $z \in X$, then $\mu'$ is quasiinvariant [Re1, Proposition 1.3.6]. We give a different proof. First of all,
$$\alpha^{\mu'} = \int \alpha^x d\mu'(x) = \mu^z \times \mu^z,$$
so $\alpha^{\mu'}$ is symmetric. Hence $\underline{\Delta}_{\mu'} = 1$. Next we consider $\lambda^{\mu'} = \int \int \beta^x_y d\mu^z(x) \, d\mu^z(y)$. Since $Gz$ is locally compact, there is a Borel function $c : [z] \to Gz$ such that for every $x \in [z]$ we have $c(x) \in xGz$. The value of $c(z)$ can be taken to be $z$. We can use $c$ to define a Borel isomorphism $\psi : G[z] \to [z] \times G(z) \times [z]$ by
$$\psi(\gamma) = (r(\gamma), c(r(\gamma))^{-1}\gamma c(s(\gamma)), s(\gamma)).$$
By the uniqueness of Haar measure, as above, we see that $\psi$ always carries $\beta^x_y$ to a positive multiple of $\epsilon^x \times \beta^z \times \epsilon^y$, and that multiple is a Borel function of the pair $(x, y)$. Hence $\psi$ carries $\lambda^{\mu'}$ to a measure equivalent to $\mu^z \times \beta^z \times \mu^z$. It follows that $\lambda^{\mu'}$ is quasisymmetric, as needed.

Since $\lambda$ is a Haar system, we know that if $K$ is a compact set in $G$ then the function $x \mapsto \lambda^x(K)$ is bounded. We will use the formula for $\lambda^x$ in terms of $\alpha^x$ to prove that $x \mapsto \alpha^x(\theta(K))$ is also bounded, and also that $\mu^x$ is finite on compact sets for the quotient topology on $[x]$. Let $F$ be the function used above to make a choice of Haar measures $\beta^y$. If $S$ is the support of $F$, then $\beta^y(S) \geq 1$ for every $y \in X$. To prove the boundedness statement above, let $K$ be a compact subset of $G$ and set $K_1 = K(s(K)S)$. Because both factors are compact, so is $K_1$, so $x \mapsto \lambda^x(K_1)$ is bounded. For $(x, y) \in \theta(K)$, choose $\gamma \in K$ such that $\theta(\gamma) = (x, y)$. Then $\gamma S \subseteq K_1$, so $\beta^x_y(K_1) \geq 1$. Hence
$$\begin{aligned}\lambda^x(K_1) &= \int \beta^x_y(K_1) \, d\alpha^x(x, y) \\ &\geq \int_{s(xK)} \beta^x_y(K_1) \, d\alpha^x(x, y) \\ &\geq \alpha^x(\theta(K)).\end{aligned}$$
For the second assertion, suppose that $C$ is a compact set in $[x]$ for the quotient topology. Since $xG$ is locally compact and $s$ is continuous and open from $xG$ to $[x]$, there is a compact set $K$ contained in $xG$ whose image contains $C$. Then $\theta(K) \subseteq xR$, so the boundedness result just proves that $\mu^x(s(K)) = \alpha^x(\theta(K))$ is finite. Hence $\mu^x$ is $\sigma$-finite.

It is also true that finiteness of $\mu^x$ on compact sets forces $\lambda^x$ to be finite on compact sets, by an argument on page 7 of [Re3].

Define $\mathcal{M}_{\theta c}(R)$ to be the space of bounded Borel functions on $R$ that vanish off sets of the form $\theta(K)$, where $K$ is a compact subset of $G$. Now we know that $\mathcal{M}_{\theta c}(R) \subseteq I(R, \lambda)$, and it is not difficult to show that $\mathcal{M}_{\theta c}(R)$ is a $*$-subalgebra of $I(G, \lambda)$. The definition of this algebra is admittedly somewhat unusual, but the algebra will serve a useful purpose in proving the main step along one way to prove the completeness of the Fourier-Stieltjes algebra of $G$. The point is that $R$ is a kind of shadow of $G$, and we need a convolution algebra on it that is a shadow of the same kind.

## 3. Positive Definite Functions

In this section, we will characterize the functions on a locally compact groupoid that are diagonal matrix entries of unitary representations as the functions that are what we call positive definite. For this to be meaningful, we need a good definition of 'positive definite'. This is more complicated than for locally compact groups because unitary representations of locally compact groupoids can be Borel functions without being continuous. Thus we make our definition using integrals, and must even identify two functions that agree $\lambda^{\mathcal{Q}}$-almost everywhere, as defined in Section 1. In Section 4, we will need to construct a positive definite function from a parametrized family of functions, each of which is positive definite on a transitive subroupoid. Thus we prove the representation theorem in that broader context. For a locally compact groupoid that has a Haar system, the notion of positive definite function can be defined in the least restrictive way as follows:

**Definition 3.1**. Let $G$ be a locally compact groupoid and let $\lambda$ be a left Haar system on $G$. Then a bounded Borel function $p$ on $G$ is called *positive definite* iff for each $x \in X$ and each $f$ in $C_c(G)$ we have

$$\text{(P)} \qquad \int\int f(\gamma_1)\bar{f}(\gamma_2)p(\gamma_2^{-1}\gamma_1)d\lambda^x(\gamma_1)d\lambda^x(\gamma_2) \geq 0.$$

The set of all such $p$'s will be denoted by $\mathcal{P}(G)$. We say that two elements of $\mathcal{P}(G)$ are *equivalent* iff they agree $\lambda^{\mathcal{Q}}$-almost everywhere.

**Remarks**. a) In Section 2, we summarized some results of J. Renault on the structure of Haar systems, and used it to construct a specific measure $\mu^x \in \mathcal{Q}$ in the unique invariant class associated with the orbit $[x]$, so that $x \sim y$ implies $\mu^x = \mu^y$. For each $x \in X$, (P) imposes a condition on the behavior of $p$ almost everywhere on the set $G|[x]$, with respect to the measure $\lambda^{\mu^x}$. Because of the invariance of the measures $\mu^z$ this measure is determined by the orbit of $x$. Any two functions that agree a.e. with respect to all such measures also agree $\lambda^{\mathcal{Q}}$-a.e. We will prove that every element of $\mathcal{P}(G)$ agrees $\lambda^{\mathcal{Q}}$-a.e. with a diagonal matrix entry of a unitary representation of $G$, and after that work always with elements of $\mathcal{P}(G)$ in that form.

b) Since $p$ is bounded, the positivity condition, (P), holds for $f$ in $L^1(\lambda^x)$ for every $x$. Furthermore, (P) can be verified using any dense subset of each $L^1(\lambda^x)$.

c) Condition (P) could be formulated using complex measures absolutely continuous relative to $\lambda^x$. (P) could be strengthened by allowing more measures, including discrete ones. Since we are able to work with the weaker definition, we will do so. In the end we will arrive at diagonal matrix entries for unitary representations, and these are as well-behaved as possible.

d) If $G$ is a Borel groupoid and $\lambda$ is a Borel Haar system, a similar definition can be made using test functions from the spaces $L^1(\lambda^x)$.

Since we intend to show that positive definite functions are essentially the same as diagonal matrix entries of unitary representations, we begin by showing that such matrix entries are in $\mathcal{P}(G)$.

**Lemma 3.2.** Let $\pi$ be a unitary representation of $G$ on a Hilbert bundle $\mathcal{H}$, and let $\xi$ be a bounded Borel section of $\mathcal{H}$. Define $p(\gamma) = (\pi(\gamma)\xi \circ s(\gamma) \mid \xi \circ r(\gamma))$ for $\gamma \in G$. Then $p \in \mathcal{P}(G)$.

**Proof.** Fix $x \in X$ and $f \in C_c(G)$. Then for $\eta \in \mathcal{H}(x)$,

$$\left| \int f(\gamma)(\pi(\gamma)\xi \circ s(\gamma) \mid \eta) \, d\lambda^x(\gamma) \right| \leq \|f\|_1 \|\xi\|_\infty \|\eta\|,$$

so there is an element $\zeta(x) \in \mathcal{H}(x)$ such that for all $\eta \in \mathcal{H}(x)$ we have

$$\int f(\gamma)(\pi(\gamma)\xi \circ s(\gamma) \mid \eta) \, d\lambda^x(\gamma) = (\zeta(x) \mid \eta).$$

Indeed, this defines a Borel section, $\zeta$, of $\mathcal{H}$. The Borel character of $\zeta$ follows from the fact that $(\pi(\gamma)\xi_1 \circ s(\gamma) \mid \eta_1 \circ r(\gamma))$ is a Borel function of $\gamma$ whenever $\xi_1$ and $\eta_1$ are Borel sections of $\mathcal{H}$. For this section $\zeta$ the integral involved in the condition (P) is equal to $(\zeta(x) \mid \zeta(x))$, which is certainly non-negative.

**Lemma 3.3.** If $p \in \mathcal{P}(G)$, the formula

$$(\text{IP}) \qquad (f \mid g)_x = \int\int f(\gamma_1) \bar{g}(\gamma_2) p(\gamma_2^{-1}\gamma_1) \, d\lambda^x(\gamma_1) \, d\lambda^x(\gamma_2)$$

defines a semi-inner product on $L^1(\lambda^x)$. Let $\mathcal{H}(x)$ denote the Hilbert space completion of the resulting inner product space. Then $\mathcal{H}$ is a Hilbert bundle over $X$. For $\gamma_1 \in G$, define $\pi(\gamma_1)$ from $L^1(\lambda^{s(\gamma_1)})$ to $L^1(\lambda^{r(\gamma_1)})$ by $(\pi(\gamma_1)f)(\gamma) = f(\gamma_1^{-1}\gamma)$. Then $\pi$ determines a unitary representation, also denoted by $\pi$, on the bundle $\mathcal{H}$.

**Proof.** The form $(f \mid g)_x$ is clearly linear in $f$ and conjugate linear in $g$. Since the vector space is complex the Hermitian symmetry follows from positive definiteness.

Let $\mathcal{N}(x) = \{f \in L^1(\lambda^x) : (f \mid f)_x = 0\}$ and set $\mathcal{F}(x) = L^1(\lambda^x)/\mathcal{N}(x)$, the corresponding inner-product space. Write $\mathcal{H}(x)$ for the completion of $\mathcal{F}(x)$. Let $|\ |_x$ be the norm (or semi-norm) arising from $(\ |\ )_x$. For $f, g \in L^1(\lambda^x)$, $|(f \mid g)_x| \leq \|p\|_\infty \|f\|_1 \|g\|_1$, so $|f|_x \leq \|p\|_\infty^{1/2} \|f\|_1$. It follows that the image of $C_c(xG)$, which is the image of $C_c(G)$, is dense in $\mathcal{H}(x)$.

Now we want to make a Borel structure on the graph of $\mathcal{H}$, denoted by $\Gamma = \Gamma_\mathcal{H} = \{(x,\xi) : x \in X, \xi \in \mathcal{H}(x)\}$. The process used is fairly standard. First, if $f \in C_c(G)$ and $x \in X$, define $\sigma(f)(x)$ to be the element of $\mathcal{H}(x)$ represented by $f|xG$. This defines a section $\sigma(f)$ of the graph of $\mathcal{H}$. We want all $\sigma(f)$'s to be Borel sections, and that tells us how to define the Borel structure. For $f \in C_c(G)$ define $\psi_f$ on $\Gamma$ by $\psi_f(x,\xi) = (\sigma(f)(x) \mid \xi)_x$. Then give $\Gamma$ the smallest Borel structure relative to which the projection to $X$ is Borel along with all the functions $\psi_f$ ($f \in C_c(G)$). It follows from the fact that $p$ is Borel and bounded that each section $\sigma(g)$ for $g \in C_c(G)$ is indeed a Borel section. Since $G$ is second countable, there is a countable set dense in $C_c(G)$. For any countable dense set of $f$'s, the $\psi_f$'s would determine the same Borel structure as $\{\psi_f : f \in C_c(G)\}$, so the latter is standard: Apply the Gram-Schmidt process in a pointwise manner to a dense sequence of sections of the form $\sigma(f)$ to get a sequence $g_1, g_2, \ldots$ of Borel functions such that

(a) $g_n|xG$ is always in $L^1(\lambda^x)$.

(b) if $f \in C_c(G)$ and $n \geq 1$, then $x \to (\sigma(f) \mid \sigma(g_n))_x$ is a Borel function.

(c) for each $x$ the non-zero elements of $\{\sigma(g_n)(x) : n \geq 1\}$ form an orthonormal basis of $\mathcal{H}(x)$.

Then it is easy to show that $\Gamma$ is isomorphic to the disjoint union of a sequence of product bundles $X_n \times \mathcal{K}_n$, where $\{X_\infty, X_1, X_2, \ldots\}$ is a Borel partition of $X$ and each $\mathcal{K}_n$ is a Hilbert space of dimension $n$. Thus $\Gamma$ is standard because each $X_n \times \mathcal{K}_n$ is standard.

If $f \in L^1(\lambda^x)$ and $\gamma_1 : x \to y$ is in $G$, define $\pi(\gamma_1)f$ by $(\pi(\gamma_1)(f))(\gamma) = f(\gamma_1^{-1}\gamma)$ for $\gamma \in yG$. Since $\lambda$ is left invariant, $\pi(\gamma_1)f \in L^1(\lambda^y)$. Notice that $\pi(\gamma_1^{-1})$ is the inverse of $\pi(\gamma_1)$. If $g$ is another element of $L^1(\lambda^x)$, then

$$(\pi(\gamma_1)f \mid \pi(\gamma_1)g)_y = \int\int f(\gamma_1^{-1}\gamma_2) \bar{g}(\gamma_1^{-1}\gamma_3) p(\gamma_3^{-1}\gamma_2) d\lambda^y(\gamma_2) d\lambda^y(\gamma_3)$$

$$= \int\int f(\gamma_2)\bar{g}(\gamma_3)p(\gamma_3^{-1}\gamma_2)d\lambda^x(\gamma_2)d\lambda^x(\gamma_3)$$

$$= (f \mid g)_x.$$

Hence $\pi(\gamma_1)$ extends to a unitary operator from $\mathcal{H}(x)$ to $\mathcal{H}(y)$, for which we use the same notation.

To work with the bundle and with the representation we need to restrict to subsets of product spaces where the various operations are defined. There are two fibered products, $\Gamma \times' \Gamma = \{(x, \xi, x, \xi') : x \in X,\ \xi, \xi' \in \mathcal{H}(x)\}$, a subset of $\Gamma \times \Gamma$, and $G \times' \Gamma \subseteq G \times \Gamma$, defined to be $\{(\gamma, x, \xi) : s(\gamma) = x,\ \xi \in \mathcal{H}(x)\}$. Let us show that $(\gamma, x, \xi) \mapsto (r(\gamma), \pi(\gamma)\xi)$ is Borel from $G \times' \Gamma$ to $\Gamma$. The composition of this map with the projection to $X$ is clearly Borel. Let $f \in C_c(G)$ and compose the map with $\psi_f$. The value of the composition at $(\gamma, x, \xi)$ is $\psi_f(r(\gamma), \pi(\gamma)\xi) = (\pi(\gamma)^{-1}(\sigma(f)(r(\gamma))) \mid \xi)_x$. This is the value of another composition:

$$G \times' \Gamma \to \Gamma \times' \Gamma \to \mathbb{C},$$

where the first function takes $(\gamma, x, \xi)$ to $(s(\gamma), \pi(\gamma^{-1})(\sigma(f)(r(\gamma))); x, \xi)$ and the second is the inner product function. The first function is Borel if each component is, so let us see that the first component is a Borel function of $\gamma$. Composition of it with projection is $s$ and hence Borel. If $g \in C_c(G)$,

$$\psi_g(s(\gamma), \pi(\gamma^{-1})(\sigma(f)(r(\gamma)))) = \int \int \bar{f}(\gamma\gamma_2) g(\gamma_1) p(\gamma_2^{-1}\gamma_1)\, d\lambda^{s(\gamma)}(\gamma_1)\, d\lambda^{s(\gamma)}(\gamma_2)$$
$$= \int_{G \times G \times G} F\, d(\epsilon^\gamma \times \lambda^{s(\gamma)} \times \lambda^{s(\gamma)}),$$

where $F$ is the function that is 0 at $(\gamma_0, \gamma_1, \gamma_2)$ unless $s(\gamma_0) = r(\gamma_1) = r(\gamma_2)$ and then its value is $\bar{f}(\gamma_0\gamma_1) g(\gamma_2) p(\gamma_2^{-1}\gamma_1)$ ($F$ is Borel). A fairly standard argument then shows that $\gamma \mapsto \psi_g(s(\gamma)), \pi(\gamma^{-1})(\sigma(f)(r(\gamma)))$ is Borel, as desired.

To show that the inner product is Borel on $\Gamma \times' \Gamma$, we use the functions $g_n$ used to show that the bundle is standard. Indeed,

$$((x, \xi) \mid (x, \eta))_x = \sum_{n \geq 1} \psi_{g_n}(x, \xi) \bar{\psi}_{g_n}(x, \eta)$$

which is a Borel function. It follows that $(\gamma, x, \xi) \mapsto \psi_f(r(\gamma), \pi(\gamma)\xi)$ is Borel, as needed.

This completes the construction of a unitary representation from a positive definite function. *From now on, subscripts will be used on inner products and norms associated with such bundles only when necessary to make clear which space is involved.* Our next task is to find a (cyclic) section $\xi_p$ such that $p(\gamma) = (\pi(\gamma)\xi_p(s(\gamma)) \mid \xi_p(r(\gamma)))$ for $\lambda^{\mathcal{Q}}$-almost every $\gamma \in G$.

The argument can be outlined as follows. Given a $\mu \in \mathcal{Q}$ we let $\mathcal{H}(\mu)$ denote the Hilbert space of square integrable sections of $\mathcal{H}$, which is sometimes written $L^2(\mu, \mathcal{H})$. There is no loss of generality in assuming that $\mu$ is a probability measure, since changing to an equivalent measure produces an equivalent representation. The representation $\pi$ of $G$ can be integrated to give a representation of $C_c(G)$ on $\mathcal{H}(\mu)$, denoted by $\pi_\mu$, using the formulation of Hahn, rather than that of Renault. The definition is given below. If $u_1, u_2, \ldots$ is a symmetric approximate unit for $C_c(G)$, the sequence of sections $\sigma(u_1), \sigma(u_2), \ldots$ has a subsequence that converges weakly to a section $\xi_\mu$ such that for $f \in C_c(G)$ we have $\pi_\mu(f)\xi_\mu = \sigma(f)$, and the matrix entry made from $\pi$ and $\xi_\mu$ agrees with $p$ a.e. relative to $\lambda^\mu$. To get a section $\xi_p$ not depending on $\mu$, we observe that if we had such a $\xi_p$, then for $f \in C_c(G)$ we would get $\int f p\, d\lambda^x = (\sigma(f) \mid \xi_p)_x$. Thus we consider the set $D(p)$ of those $x \in X$ for which $\int f p\, d\lambda^x$, as a linear function of $f$ in $C_c(G)$, 'extends' to a bounded linear functional on $\mathcal{H}(x)$. We need to know that $D(p)$ is conull for every $\mu \in \mathcal{Q}$, and this follows from the existence of $\xi_\mu$. We let $\xi_p(x)$ be the vector representing that linear functional, and verify that $\xi_p$ is the section we wanted.

Before giving details, we introduce the space $L^{1,2}(\lambda, \mu)$, consisting of those Borel functions $f$ for which

$$\|f\|_{1,2}^2 = \int \left( \int |f(\gamma)|\, d\lambda^x(\gamma) \right)^2 d\mu(x) < \infty.$$

Now, begin by taking $\mathcal{H}(\mu)$ as defined above, and observe that for $f \in L^{1,2}(\lambda, \mu)$, the section $\sigma(f)$ is in $\mathcal{H}(\mu)$, and $\|\sigma(f)\| \leq \|p\|_\infty^{1/2} \|f\|_{1,2}$, so that $f_n \to f$ in $L^{1,2}(\lambda, \mu)$ implies $\sigma(f_n) \to \sigma(f)$ in $\mathcal{H}(\mu)$. To make

the proof work, we must integrate the representation $\pi$ to get $\pi_\mu$ having the property that for $f, g \in C_c(G)$, $\pi_\mu(f)(\sigma(g)) = \sigma(f*g)$. This can be done if we use the method of [H2], which applies to $I(G, \lambda)$, which is a subspace of $L^{1,2}(\lambda, \mu)$ because $\mu$ is a probability measure.

For $f \in I(G, \lambda)$ we define $\pi_\mu(f)$ by saying that for sections $\xi, \eta$ in $\mathcal{H}(\mu)$ we have

$$(\pi_\mu(f)\xi \mid \eta) = \int f(\gamma)(\pi(\gamma)\xi(s(\gamma)) \mid \eta(r(\gamma)))\, d\lambda^\mu(\gamma).$$

The integral defines a bounded sesquilinear form, so the formula defines a bounded operator $\pi_\mu(f)$. It is proved in [H2] that $\pi_\mu$ is a bounded representation of $I(G, \lambda)$. If $f \in I(G, \lambda)$ and $\xi \in \mathcal{H}(\mu)$, then $\pi_\mu(f)\xi$ is represented by a section whose value at almost every $x$ is $\int f(\gamma)\pi(\gamma)\xi(s(\gamma))\, d\lambda^x(\gamma)$, where the integral is defined weakly. If $g \in I(G, \lambda)$, we have

$$\begin{aligned}(\pi_\mu(f)\sigma(g) \mid \eta) &= \int f(\gamma)(\pi(\gamma)(\sigma(g)(s(\gamma))) \mid \eta(r(\gamma)))d\lambda^\mu(\gamma) \\ &= \int\int f(\gamma)(\pi(\gamma)(\sigma(g)(s(\gamma)))\mid \eta(x))d\lambda^x(\gamma)\, d\mu(x) \\ &= \int (\sigma(f*g)(r(\gamma)) \mid \eta(r(\gamma))d\lambda^\mu(\gamma),\end{aligned}$$

because $\pi(\gamma)(\sigma(g)(s(\gamma)))$ is represented by a function on $r(\gamma)G$ whose value at a point $\gamma_1$ is $g(\gamma^{-1}\gamma_1)$.

**Lemma 3.4.** Let $G$ be a locally compact groupoid with a Haar system $\lambda$. Suppose that $\mu \in \mathcal{Q}$ is a probability measure. If $p$ is a positive definite function on $G$ and $(\mathcal{H}, \pi)$ is constructed from $p$ as in the proof of Lemma 3.3, then there is a section $\xi_\mu \in \mathcal{H}(\mu)$ such that

(a) $|\xi_\mu(x)|^2_x \leq \|p\|_\infty$ for $x \in X$

(b) $\pi_\mu(f)\xi_\mu = \sigma(f)$ for $f \in C_c(G)$

(c) $p(\gamma) = (\pi(\gamma)\xi_\mu(s(\gamma)) \mid \xi_\mu(r(\gamma)))$ a.e. $d\lambda^\mu(\gamma)$.

**Proof** Let $u_1, u_2, \ldots$ be a symmetric approximate unit for $G$. Then $|\sigma(u_i)(x)| \leq \|p\|_\infty^{1/2}$ for each $x$ and $i$, so $\|\sigma(u_i)\| \leq \|p\|_\infty^{1/2}$ for each $i$. Thus $\sigma(u_1), \sigma(u_2), \ldots$ has a subsequence converging weakly to a vector $\xi_\mu \in \mathcal{H}(\mu)$. We may suppose that subsequence is $\sigma(u_1), \sigma(u_2), \ldots$. If, for every Borel set $E$ in $X$, $P(E)$ is the projection of $\mathcal{H}(\mu)$ onto the subspace determined by sections that vanish off $E$, then $P(E)\sigma(u_n)$ converges weakly to $P(E)\xi_\mu$, which has norm at most $(\|p\|_\infty \mu(E))^{1/2}$. It follows that $|\xi_\mu(x)|^2_x \leq \|p\|_\infty$ for a.e. x, and we can change $\xi_\mu$ to make it true for all $x$. For $f \in C_c(G)$, $f*u_i \to f$ uniformly and all these functions vanish off a fixed compact set. Thus $f*u_i \to f$ in $L^{1,2}$, and $\sigma(f*u_i) \to \sigma(f)$ in $\mathcal{H}(\mu)$. Hence $\pi_\mu(f)\sigma(u_i)$ converges to $\sigma(f)$. We also know that $\pi_\mu(f)$ is a bounded operator, so $\pi_\mu(f)\sigma(u_i)$ converges weakly to $\pi_\mu(f)\xi_\mu$. Hence $\pi_\mu(f)\xi_\mu = \sigma(f)$, as elements of $\mathcal{H}(\mu)$.

It follows from this that if $f, g \in C_c(G)$, then $(\sigma(f) \mid \sigma(g)) = (\pi_\mu(f)\xi_\mu \mid \pi_\mu(g)\xi_\mu)$ and this can be written as

$$\int\int\int (\pi(\gamma_1)\xi_\mu(s(\gamma_1)) \mid \pi(\gamma_2)\xi_\mu(s(\gamma_2)))f(\gamma_1)\bar{g}(\gamma_2)d\lambda^x(\gamma_1)d\lambda^x(\gamma_2)d\mu(x)$$

which is equal to

$$\int\int\int (\pi(\gamma_2^{-1}\gamma_1)\xi_\mu(s(\gamma_1)) \mid \xi_\mu(s(\gamma_2)))f(\gamma_1)\bar{g}(\gamma_2)d\lambda^x(\gamma_1)d\lambda^x(\gamma_2)d\mu(x).$$

If $h \in C_c(X)$ we can replace $f$ in this calculation by $hf = (h \circ r)f$. From this it follows that if $f, g \in C_c(G)$ then for $\mu$-almost every $x$ we have

$$(\sigma(f) \mid \sigma(g))_x = \int\int (\pi(\gamma_2^{-1}\gamma_1)\xi_\mu(s(\gamma_1)) \mid \xi_\mu(s(\gamma_2)))f(\gamma_1)\bar{g}(\gamma_2)d\lambda^x(\gamma_1)d\lambda^x(\gamma_2).$$

By the definition of $(\mid)_x$, this shows that for $\mu$-almost every $x$,

(∗) $$p(\gamma_2^{-1}\gamma_1) = (\pi(\gamma_2^{-1}\gamma_1)\xi_\mu(s(\gamma_1)) \mid \xi_\mu(s(\gamma_2)))$$

is true for $\lambda^x \times \lambda^x$-almost all pairs $(\gamma_1, \gamma_2)$. For each such $x$, for $\lambda^x$-almost every $\gamma_2$, the formula $(*)$ is true for $\lambda^x$-almost every $\gamma_1$, i.e., $p(\gamma) = (\pi(\gamma)\xi_\mu(s(\gamma)) \mid \xi_\mu(r(\gamma)))$ for $\lambda^{s(\gamma_2)}$–almost all $\gamma$. Indeed, the set $\{s(\gamma_2) : (*) \text{ holds for } \lambda^x\text{-almost every } \gamma_1\}$ is conull in $[x]$.

**Theorem 3.5.** Let $G$ be a locally compact groupoid and let $\lambda$ be a Haar system on $G$. If $p$ is a positive definite function on $G$ and $(\mathcal{H}, \pi)$ is the associated unitary $G$-bundle over $X$, then there is a bounded section $\xi_p$ of $\mathcal{H}$ such that if $\mu \in \mathcal{Q}$, then

(a) $p(\gamma) = (\pi(\gamma)\xi_p(s(\gamma)) \mid \xi_p(r(\gamma)))$ a.e. $d\lambda^\mu(\gamma)$

(b) if $f \in C_c(G)$, then $\pi_\mu(f)\xi_p = \sigma(f)$ in $\mathcal{H}(\mu)$.

If $p$ is continuous, then $\xi_p$ can be chosen to be continuous and $p(\gamma) = (\pi(\gamma)\xi_p(s(\gamma)) \mid \xi_p(r(\gamma)))$ for all $\gamma$.

**Proof** Define $D = D(p) = \{x \in X : f \mapsto \int fp d\lambda^x = \lambda^x(fp) \text{ extends from } C_c(G) \text{ to give a bounded linear functional on } \mathcal{H}(x) \text{ of norm at most } \|p\|_\infty^{1/2}\}$. For each $f \in C_c(G)$, $\lambda(fp)$ and $x \mapsto (f \mid f)_x$ are Borel functions, and boundedness can be tested on a countable dense set, so $D$ is a Borel set. For $x \in D$, there is a unique vector $\xi_p(x) \in \mathcal{H}(x)$ such that $(\sigma(f)(x) \mid \xi_p(x))_x = \lambda^x(fp)$ for $f \in C_c(G)$, and if we let $\xi_p(x) = 0$ for $x \notin D$, $\xi_p$ is a Borel section of $\mathcal{H}$, bounded by $\|p\|_\infty^{1/2}$. We need to show that $D$ is $\mathcal{Q}$-conull, i.e., conull for each $\mu \in \mathcal{Q}$.

Let $\mu \in \mathcal{Q}$. Then there is a $\xi_\mu \in \mathcal{H}(\mu)$ satisfying (a), (b), (c) of Lemma 3.4 and thus for each $f \in C_c(G)$ we have
$$(\sigma(f)(x) \mid \xi_\mu(x))_x = ((\pi_\mu(f)\xi_\mu(x) \mid \xi_\mu(x))_x = \lambda^x(fp)$$
for $\mu$-a.e. $x$. Since bounded linear functionals are determined by their values on a countable dense set, and since boundedness of a linear functional can be tested on a countable dense set, there is a $\mu$-conull set $D_\mu$ such that for $x \in D_\mu$ and $f \in C_c(G)$,
$$(\sigma(f)(x) \mid \xi_\mu(x))_x = \lambda^x(fp).$$

Thus $D_\mu \subseteq D$, from which it follows that $D$ is $\mu$-conull and $\xi_p(x) = \xi_\mu(x)$ $\mu$-a.e. This fact and Lemma 3.3 combine to establish the truth of statements a) and b) in the theorem. By the definitions of $D$ and $\xi_p$, it follows that $\xi_p$ is bounded by $\|p\|_\infty^{1/2}$.

**Remark**. This shows that we can replace $p(\gamma)$ by $(\pi(\gamma)\xi_p(s(\gamma)) \mid \xi_p(r(\gamma))$, i.e., assume $p$ has that form, as for groups.

To complete the proof, we show first that if $p$ is continuous, then $D = X$. Again take a $\mu \in \mathcal{Q}$ and the section $\xi_\mu$. We have $\lambda^x(fp) = (\sigma(f)(x) \mid \xi_\mu(x))_x$ for $\mu$-a.e. $x$, and for such $x$'s,
$$|\lambda^x(fp)| \leq |\sigma(f)(x)|_x \|\xi_\mu(x)\|$$
$$\leq |\sigma(f)(x)|_x \|p\|_\infty^{1/2}.$$

Since $p$ is continuous, both $\lambda(fp)$ and $x \mapsto (\sigma(f), \sigma(f))_x$ are continuous, so this estimate holds on the support of $\mu$. The supports of the $\mu$'s in $\mathcal{Q}$ fill $X$, so
$$|\lambda^x(fp)| \leq |\sigma(f)(x)| \|p\|_\infty^{1/2}$$
for all $f$ in $C_c(G)$ and all $x$. Thus $D = X$. Now $p(\gamma) = (\pi(\gamma)\xi_p(s(\gamma)) \mid \xi_p(r(\gamma)))$ a.e. $d\lambda^\mu(\gamma)$ for every $\mu$, so it will end the proof if we can show that the second function is continuous. By a partition of unity argument, this will follow if we can show that $(\pi(\gamma)\xi(s(\gamma)) \mid \xi(r(\gamma)))$ is a continuous function of $\gamma$ for every continuous section $\xi$ of compact support. In fact we can reduce to considering $\xi = \sigma(f)$ for $f \in C_c(G)$, by using partitions of unity and uniform limits. Then we have
$$(\pi(\gamma)\xi(s(\gamma)) \mid \xi(r(\gamma))) = \int\int f(\gamma^{-1}\gamma_1)\bar{f}(\gamma_2)p(\gamma_2^{-1}\gamma_1)d\lambda^{r(\gamma)}(\gamma_1)d\lambda^{r(\gamma)}(\gamma_2).$$

Continuity of this function of $\gamma$ can be deduced by applying the following easy lemma and a variant of it using the second coordinate projection, because the integrands can be extended to functions satisfying the hypotheses of the lemma.

**Lemma 3.6.** Suppose $G$ is a locally compact groupoid with a Haar system $\lambda$ and let $F$ be a continuous complex valued function on $G \times G$. Let $p_1 : G \times G \to G$ be the first coordinate projection. Suppose that for every compact set $C \subseteq G$ the set $p_1^{-1}(C) \cap \mathrm{supp}(F)$ is compact. Then, $\int F(\gamma, \gamma_2) d\lambda^{r(\gamma)}(\gamma_2)$ is a continuous function of $\gamma$.

We have an existence theorem, but we should show that the results are essentially the same for any two equivalent elements of $\mathcal{P}(G)$.

**Theorem 3.7.** Suppose that $p, q \in \mathcal{P}(G)$ and that $p = q$ $\lambda^{\mathcal{Q}}$-a.e. Then the associated representations $(\mathcal{H}_p, \pi_p)$ and $(\mathcal{H}_q, \pi_q)$ are the same, and the sections $\xi_p$ and $\xi_q$ agree $\mathcal{Q}$-a.e.

**Proof.** Let $z \in X$ and consider the inner products on $L^1(\lambda^z)$ defined using $p$ and $q$. Denote them by $(\ |\ )_p$ and $(\ |\ )_q$ respectively. To prove they are the same, it will suffice to show that $p(\gamma_2^{-1}\gamma_1) = q(\gamma_2^{-1}\gamma_1)$ for $\lambda^z \times \lambda^z$-almost every pair $(\gamma_1, \gamma_2)$, because the inner products are defined by double integrals using these functions and measures.

Let $\mu_1^z$ be a quasiinvariant probability measure equivalent to the measure $\mu^z$ that was associated with the orbit $[z]$ near the end of Section 2. Let $E$ be the set of $x \in X$ for which $p = q$ a.e. relative to $\lambda^x$. Then $\mu^z(E) = 1$, so $\{\gamma : s(\gamma) \in E\} = GE$ is $\lambda^z$-conull. If $\gamma \in zGE$, $\{\gamma_1 \in zG : p(\gamma^{-1}\gamma_1) = q(\gamma^{-1}\gamma_1)\}$ is conull relative to $\lambda^z$ by translation invariance of the Haar system. By the Fubini Theorem, we get the desired agreement a.e.

This shows that the Hilbert bundles $\mathcal{H}_p$ and $\mathcal{H}_q$ are identical, and since the formula for the representation is just left translation in each case, the representations are the same.

To show that the sections $\xi_p$ and $\xi_q$ agree $\mathcal{Q}$-a.e., we resort to the definitions, namely, $\xi_p(x)$ and $\xi_q(x)$ are determined by the fact that for $f \in C_c(G)$

$$(\sigma(f) \mid \xi_p(x)) = \lambda^x(fp)$$

and

$$(\sigma(f) \mid \xi_q(x)) = \lambda^x(fq).$$

Let $F$ be the set of $x \in X$ for which $\xi_p(x) = \xi_q(x)$. Since the two sections are Borel, $F$ is a Borel set. We need to show that if $\mu \in \mathcal{Q}$, then $\mu(X \backslash F) = 0$. We know that for each $f \in C_c(G)$ the two functions $\lambda(fp)$ and $\lambda(fq)$ agree almost everywhere relative to $\mu$. Let $\mathcal{C}$ be a countable dense set in $C_c(G)$, and let $N$ be a $\mu$-null set such that $x \notin N$ and $f \in \mathcal{C}$ imply $\lambda^x(fp) = \lambda^x(fq)$. Since $p$ and $q$ are bounded, this equality is preserved under limits in $C_c(G)$, so it holds for $x \notin N$ and all $f \in C_c(G)$. Thus $F$ contains the complement of $N$, as desired.

**Theorem 3.8.** Sums and products of positive definite functions are positive definite.

**Proof** This is immediate from the existence of direct sums and tensor products of representations, because of Theorem 3.5.

Now let us consider the enlarging the space from which we construct the fibers of the Hilbert bundle $\mathcal{H}_p$ using the positive definite function $p$. In later sections it will be convenient to replace the algebra $C_c(G)$ by the larger algebra $C_c(G, \bar{X})$, an algebra of kernels introduced in Section 1, and we will need to know that using the latter in our construction does not change the fibers in that bundle.

**Definition 3.9.** Let $G$ be a locally compact groupoid and let $\lambda$ be a left Haar system on $G$. Then a bounded Borel function $p$ on $G$ is called *strictly positive definite* if for each $x \in X$ and each $\nu$ in $C_c(G, \bar{X})$ we have

(P')
$$\int \int p(\gamma_2^{-1}\gamma_1) d\nu^x(\gamma_1) d\bar{\nu}^x(\gamma_2) \geq 0.$$

The set of all such $p$'s will be denoted by $\mathcal{P}'(G)$. Two functions $p, q \in \mathcal{P}'(G)$ will be called *equivalent* iff they agree $\lambda^{\mathcal{Q}}$-almost everywhere on $G$ and their restrictions to $X$ agree $\mathcal{Q}$-almost everywhere.

**Remark.** Equivalence of two elements of $\mathcal{P}'(G)$ just guarantees that they lead to the same semi-inner products on $C_c(G, \bar{X})$. Part of the proof of this is given by Theorem 3.7, and the other follows from the fact that the rest of inner product involves integration over $X$ at least once.

Of course we have $\mathcal{P}'(G) \subseteq \mathcal{P}(G)$, and would like to know that the sets are equal. Strictly speaking, this is not true because a function $p$ can satisfy condition (P) and be negative everywhere on $X$ unless there is a $\mu \in \mathcal{Q}$ such that $\lambda^\mu(X) > 0$. Actions by non-discrete groups give rise to groupoids for which $\mathcal{Q}$ contains no such $\mu$. However, we have proved that every equivalence class in $\mathcal{P}(G)$ contains a diagonal matrix entry. Thus a kind of reverse of the containment would follow from the analog of Lemma 3.2, namely Lemma 3.11 below showing that diagonal matrix entries are in $\mathcal{P}'(G)$. This meaning of the reverse containment would be that every class in $\mathcal{P}(G)$ contains an element of $\mathcal{P}'(G)$, or that diagonal matrix entries are in $\mathcal{P}'(G)$.

However, there is another natural question that also should be answered. If two diagonal matrix entries are equivalent in $\mathcal{P}(G)$, are they equivalent in $\mathcal{P}'(G)$? The affirmative answer is given in Lemma 3.13.

**Lemma 3.10.** Let $\pi$ be a unitary representation of $G$ on the Hilbert bundle $\mathcal{K}$, and let $\xi$ be a bounded Borel section of $\mathcal{K}$. Define $p(\gamma) = (\pi(\gamma)\xi \circ s(\gamma) \mid \xi \circ r(\gamma))$ for $\gamma \in G$. Then $p \in \mathcal{P}'(G)$.

**Proof.** As in the proof of Lemma 3.2, there is a section, $\zeta$, of the bundle such that for each $x \in X$ and every $\eta \in \mathcal{K}(x)$ we have

$$\int f(\gamma)(\pi(\gamma)\xi \circ s(\gamma) \mid \eta)\, d\lambda^x(\gamma) = (\zeta(x) \mid \eta).$$

If $g \in C(\bar{X})$, $x \in X$, and $\eta \in \mathcal{K}(x)$, then

$$\int g(\gamma)(\pi(\gamma)\xi \circ s(\gamma) \mid \eta)\, d\epsilon^x(\gamma) = g(x)(\xi(x) \mid \eta).$$

If $\nu = f\lambda + g\epsilon$, these show that the integral involved in the condition (P') is equal to $(\zeta(x) + g(x)\xi(x) \mid \zeta(x) + g(x)\xi(x))$, which is certainly non-negative.

**Corollary 3.11.** Every equivalence class in $\mathcal{P}(G)$ contains an element of $\mathcal{P}'(G)$.

**Lemma 3.12.** Let $(\mathcal{H}, \pi)$ be a representation of $G$, let $u_1, u_2, \ldots$ be a symmetric approximate unit [M-R-W], as described in Theorem 1.1, and let $\xi$ be a bounded Borel section of $\mathcal{H}$. Suppose that $\mu \in \mathcal{Q}$, and let $\pi_\mu$ be the integrated form of $\pi$ as defined just before the statement of Lemma 3.3. Then $\pi_\mu(u_n)\xi \to \xi$ as $n \to \infty$.

**Proof.** By construction of the functions $u_n$, $\|u_n\|_I \leq 1$ for each $n$, so every $\pi_\mu(u_n)$ has norm at most 1. Hence it suffices to find a dense set of vectors satisfying the conclusion. Each vector of the form $\pi_\mu(g)\eta$ satisfies the conclusion, and hence vectors in the linear span of the set of such vectors do also. That linear span is dense.

**Lemma 3.13.** Suppose that $\pi$ and $\pi_1$ are representations of $G$ on bundles $\mathcal{K}$ and $\mathcal{K}_1$, and that $\xi$ and $\xi_1$ are bounded Borel sections of $\mathcal{K}$ and $\mathcal{K}_1$. If $(\pi(\gamma)\xi \circ s(\gamma) \mid \xi \circ r(\gamma)) = (\pi_1(\gamma)\xi_1 \circ s(\gamma) \mid \xi_1 \circ r(\gamma))$ for almost every $\gamma$ relative to $\lambda^\mathcal{Q}$, then $(\xi(x) \mid \xi(x)) = (\xi_1(x) \mid \xi_1(x))$ for almost every $x \in X$ relative to $\mathcal{Q}$.

**Proof.** Since $\xi$ and $\xi_1$ are Borel sections the set $E$ of $x \in X$ for which $(\xi(x) \mid \xi(x)) = (\xi_1(x) \mid \xi_1(x))$ is a Borel set. We need to prove that for $\mu \in \mathcal{Q}$, $\mu(E) = 1$. The hypothesis implies that for $f \in \mathcal{M}_c(G)$ we have $(\pi_\mu(f)\xi \mid \xi) = (\pi_{1,\mu}(f)\xi_1 \mid \xi_1)$, these being inner products associated with the integrated forms using Hahn's method (cf. Lemma 3.4, and the paragraph before it). Let $\varphi$ and $\varphi_1$ be the representations of $C(\bar{X})$ by multiplication on the sections of $\mathcal{K}$ and $\mathcal{K}_1$. Then it follows from the discussion following the statement of Theorem 1.2 that for $h \in C(\bar{X})$ and $f \in C_c(G)$

$$(\varphi(h)\pi_\mu(f)\xi \mid \xi) = (\varphi_1(h)\pi_{1,\mu}(f)\xi_1 \mid \xi_1).$$

Now for the $f \in C_c(G)$ take the terms of a symmetric approximate unit, to see that for all $h \in C(\bar{X})$,

$$(\varphi(h)\xi \mid \xi) = (\varphi_1(h)\xi_1 \mid \xi_1).$$

This means that for all $h \in C(\bar{X})$,

$$\int h(x)(\xi(x) \mid \xi(x))\, d\mu(x) = \int h(x)(\xi_1(x) \mid \xi_1(x))\, d\mu(x).$$

Thus $E$ is indeed $\mu$-conull.

After this, we will always take elements of $\mathcal{P}(G)$ or $\mathcal{P}'(G)$ to be diagonal matrix entries, and understand that they are determined a.e. on $X$ as well as on $G$.

## 4. Complete Positivity

In this section we introduce second and third ways to view elements of $\mathcal{P}(G)$, namely in terms of completely positive mappings. Theorem 4.1 is a first step toward getting Banach algebras of completely bounded maps on $M^*(G)$ and on $C^*(G)$. In Section 1 we obtained $C^*(G)$ by completing $C_c(G)$, and defined $\omega$ to be the direct sum of the (cyclic) representations of $C^*(G)$ that arise from normalized positive linear functionals on $C^*(G)$. Let $\mathcal{H}_\omega$ be the Hilbert space of $\omega$. By a theorem of Renault, stated in Section 1, each representation of $C_c(G)$ can be gotten by integrating a unitary representation of $G$. Thus $\omega|C_c(G)$ is also a direct sum of certain representations $\pi^\mu$. The process of integration allows us to regard each $\pi^\mu$, and hence $\omega$, as a representation of either $\mathcal{M}_c(G)$ or $C_c(G)$. We call $\omega$ the universal representation of $G$. We also defined $M^*(G)$ to be the operator norm closure of $\omega(\mathcal{M}_c(G))$, and notice that $C^*(G)$ is isomorphic to the norm closure of $\omega(C_c(G))$. If $G$ is a group, of course $M^*(G) = C^*(G)$, but these two algebras can be different for groupoids.

**Theorem 4.1** Let $p$ be a positive definite function on $G$. Let $\omega$ be the universal representation of $G$, and define $T_p(\omega(f)) = \omega(pf)$ for $f \in \mathcal{M}_c(G)$. Then $T_p$ extends to a completely positive map of $M^*(G)$ to $M^*(G)$ with completely bounded norm equal to the $\mathcal{Q}$-essential supremum of $\{p(x) : x \in X\}$.

**Proof.** (We modify a proof from [W4] for groupoids.) We remind the reader that this $\mathcal{Q}$-essential supremum is the infimum of $\{B : \text{if } \mu \in \mathcal{Q}, \text{ then } p \le B \ \mu\text{-a.e.}\}$. Also, in working with $\omega$ we will use its construction as a direct sum.

We will need to find a formula for $T_p$, in order to prove that the mapping is completely positive. For this we begin with two vectors $\xi, \eta$ in one summand of $\mathcal{H}_\omega$ given by an integrated representation $\pi^\mu$. This means that we begin with a measure $\mu \in \mathcal{Q}$ and a Hilbert bundle $\mathcal{K}$ over $X$. The subspace of $\mathcal{H}_\omega$ in question is $L^2(\mu; \mathcal{K})$, and the restriction of $\omega$ to this subspace is the integrated form of a representation, $\pi$, of $G$. We are using Renault's form here, as described in Section 1: take $\nu = \int \lambda^x d\mu(x)$ and $\nu_0 = \Delta_\mu^{-1/2} \nu$. Then for $f \in \mathcal{M}_c(G)$,

$$(T_p \omega(f) \xi \mid \eta) = (\omega(pf) \xi \mid \eta)$$
$$= \int p(\gamma) f(\gamma) (\pi(\gamma) \xi(s(\gamma)) \mid \eta(r(\gamma))) d\nu_0(\gamma)$$
$$= (\pi^\mu(pf) \xi \mid \eta).$$

By Theorem 3.4 there are a Hilbert bundle $\mathcal{K}_p$ on $X$, a (unitary) Borel representation $\pi_p$ of $G$ on $\mathcal{K}_p$ and a bounded Borel section $\xi_p$ of $\mathcal{K}_p$ such that $p(\gamma) = (\pi_p(\gamma) \xi_p \circ s(\gamma) \mid \xi_p \circ r(\gamma))$ for $\lambda^\mu$-a.e. $\gamma \in G$. By Theorem 3.7, $\pi_p$ is unique, and the section $\xi_p$ is determined $\mathcal{Q}$-a.e. Thus we can continue the calculation from above as follows:

$$= \int f(\gamma) (\pi_p(\gamma) \xi_p \circ s(\gamma) \mid \xi_p \circ r(\gamma))(\pi(\gamma) \xi \circ s(\gamma) \mid \eta \circ r(\gamma)) \, d\nu_0(\gamma)$$
$$= \int f(\gamma) ((\pi_p \otimes \pi(\gamma))(\xi_p \otimes \xi) \circ s(\gamma) \mid (\xi_p \otimes \eta) \circ r(\gamma)) d\nu_0(\gamma)$$
$$= ((\pi_p \otimes \pi)(f)(\xi_p \otimes \xi) \mid \xi_p \otimes \eta)$$
$$= ((\pi_p \otimes \omega)(f)(\xi_p \otimes \xi) \mid \xi_p \otimes \eta).$$

Here $\xi_p \otimes \xi$ and $\xi_p \otimes \eta$ are in $L^2(\mu; \mathcal{K}_p \otimes \mathcal{K})$. In summary we have

$$(T_p \omega(f) \xi \mid \eta) = ((\pi_p \otimes \omega)(f) V_{p, \mu, \mathcal{K}} \xi \mid V_{p, \mu, \mathcal{K}} \eta)$$

where $V_{p, \mu, \mathcal{K}} : L^2(\mu; \mathcal{K}) \to L^2(\mu; \mathcal{K}_p \otimes \mathcal{K})$ is defined by $V_{p, \mu, \mathcal{K}} \xi = \xi_p \otimes \xi$. This is a bounded operator because the section $\xi_p$ is bounded and the usual techniques for multiplication operators apply. If we let $V_p$ be the direct sum of the operators $V_{p,\mu,\mathcal{K}}$ over all pairs $(\mu,\mathcal{K})$, we have $T_p\omega(f) = V_p^*(\pi_p \otimes \omega)(f)V_p$. A theorem of Stinespring, [P, St], shows that $T_p$ is completely positive with completely bounded norm equal to $\|V_p\|^2$. But $V_p$ is given by a tensor multiplication which behaves like a scalar multiplication operator, so

$$\|V_p\|^2 = \operatorname{ess\,sup}\{\|\xi_p(x)\|^2 : x \in X\} = \operatorname{ess\,sup}\{p(x) : x \in X\}.$$

The proof of Theorem 4.1 also proves this:

**Theorem 4.2** Let $p$ be a positive definite function on $G$, let $\mu \in \mathcal{Q}$ and let $\pi$ be a representation of $G$. Define $T_p'(\pi^\mu(f)) = \pi^\mu(pf)$ for $f \in \mathcal{M}_c(G)$. Then $T_p'$ extends to a completely positive map of the norm closure of $\pi^\mu(\mathcal{M}_c(G))$ to itself, this being the quotient of the $T_p$ defined in Theorem 4.1. The completely bounded norm of $T_p$ as an operator on $\operatorname{cl}(\pi^\mu(\mathcal{M}_c(G)))$ is the $\mu$-essential supremum of $\{p(x) : x \in X\}$.

Although the norm on the Fourier-Stieltjes algebra of a groupoid comes from its representation by completely bounded maps rather than as the Banach space dual of the $C^*$-algebra as it does for groups, the latter fact has a remnant. Here we prove just one lemma regarding that remnant.

**Lemma 4.3** Let $p$ be a positive definite function on $G$, and let $\mu$ be a probability measure in $\mathcal{Q}$. Define $\psi_{p,\mu}(\omega(f)) = \int f(\gamma)p(\gamma)d\nu(\gamma)$ for $f \in C_c(G)$, where $\nu = \int \lambda^x d\mu(x)$ (Section 1). Then $\psi_{p,\mu}$ extends to a positive linear functional on $C^*(G)$ whose norm is at most the $\mathcal{Q}$-essential supremum of $p$.

**Proof** From the definition of $\pi_\mu$ in Section 3, it follows that the integral in question is equal to $(\pi_\mu(f)\xi \mid \xi)$, where $\pi$ is the unitary representation of $G$ derived from $p$ and $\xi$ is the associated section of the Hilbert bundle. Thus this linear functional is clearly positive, and its norm is at most $\|\xi\|^2$, the square of the norm of $\xi$ in $\mathcal{H}(\mu)$, but this is at most $\|\xi\|_\infty^2$ which is the $\mathcal{Q}$-essential supremum of $p$.

Next we present a third way to think about $\mathcal{P}(G)$. It depends on using the decomposition described in Section 2 of the Haar system of $G$ over the equivalence relation $R$ associated to $G$. This decomposition is relative to the mapping $\theta = (r,s)$ of $G$ onto $R$. Since $G$ is $\sigma$-compact it follows that $R$ is $\sigma$-compact in the quotient topology. The decomposition of the Haar system involves two families of measures. First of all there is a measure $\beta_y^x$ concentrated on $xGy$ for every pair $(x,y)$ in $R$, such that each $\beta_y^y$ is a Haar measure on $yGy$ and $\beta_y^x$ is a translate of $\beta_y^y$. Then there is a Borel Haar system $\alpha$ for $R$ so that for every $x \in X$ we have

$$\lambda^x = \int \beta_y^z \, d\alpha^x(z,y).$$

There is a Borel homomorphism $\delta$ from $G$ to the positive reals such that for every $\mu \in \mathcal{Q}$ the modular homomorphisms $\Delta_\mu$ for $G$ and $\tilde{\Delta}_\mu$ for $R$ satisfy $\Delta_\mu = \delta \tilde{\Delta}_\mu \circ \theta$. For each $x \in X$ let $\mu^x$ be the measure on $X$ so that $\alpha^x = \epsilon^x \times \mu^x$. Then $x \sim y$ implies $\mu^x = \mu^y$. Thus $\alpha^{\mu^x} = \mu^x \times \mu^x$, so $\tilde{\Delta}_{\mu^x} = 1$.

Let $\mathcal{M}_{\theta c}(R)$ be the space of bounded Borel functions on $R$ supported on images under $\theta$ of compact subsets of $G$. Then $\mathcal{M}_{\theta c}(R)$ is a $*$-algebra under convolution, using the Borel Haar system $\alpha$ (see Section 2). We also extend this algebra to include $\mathcal{M}(X)$, as done in Section 1 for $\mathcal{M}_c(G)$ and $\mathcal{M}(X)$, obtaining $\mathcal{M}_{\theta c}(R,X)$ in this case.

If $\mu$ is a quasi-invariant measure on $X$, i.e., $\mu \in \mathcal{Q}$, earlier we introduced the notation $\lambda^\mu$ for $\int \lambda^x \, d\mu(x)$ and we define $\alpha^\mu$ similarly. Now we want to shorten the notation, so we write $\nu = \lambda^\mu$, $\tilde{\nu} = \alpha^\mu$, $\Delta = \Delta_\mu$, and $\tilde{\Delta} = \tilde{\Delta}_\mu$.

To integrate a unitary representation of $G$ relative to $\mu$ to make a $*$-representation of $\mathcal{M}_c(G,X)$, we use the measure $\nu_0 = \Delta^{-1/2}\nu$ and to integrate a representation of $R$ we use the measure $\tilde{\nu}_0 = \tilde{\Delta}^{-1/2}\tilde{\nu}$. For example, in the first case we have

$$(\pi^\mu(f)\xi \mid \eta) = \int f(\gamma)(\pi(\gamma)\xi \circ r(\gamma) \mid \eta \circ s(\gamma))\, d\nu_0(\gamma)$$

whenever $f \in \mathcal{M}_c(G)$ and $\xi,\eta$ are $L^2$ sections of the bundle on which $\pi$ represents $G$. This is the formulation of [Re1]. From what we have above, it follows that $\nu_0 = \int \delta^{-1/2}\beta_y^x \, d\tilde{\nu}_0(x,y)$, so there is a convenient relationship between the two measures.

For each unitary representation $\pi$ of $R$, and each $\mu \in \mathcal{Q}(R)$, we can ask whether the representation $\pi^\mu$ is cyclic, and we can define $\tilde{\omega}$ to be a direct sum formed using for summands one representative from each equivalence class of a cyclic $\pi^\mu$. Then we can write $M^*(R)$ for the norm closure of $\tilde{\omega}(\mathcal{M}_{\theta c}(R))$. These $\pi^\mu$'s extend to $\mathcal{M}_{\theta c}(R, X)$, so $\tilde{\omega}$ does also, and we let $M^*(R, X)$ be the norm closure of $\tilde{\omega}(\mathcal{M}_{\theta c}(R, X))$. As stated before, the algebra $M^*(R, X)$ is present only for its utility in proving results about $G$, and the slightly strange definition is just suited to that purpose.

If $p \in \mathcal{P}(G)$, we define a pairing of $p$ with an element $f \in \mathcal{M}_c(G)$ to give a function on $R$ as follows:

$$\langle f, p \rangle(x, y) = \int f p \delta^{-1/2} \, d\beta_y^x.$$

Since $p$ and $\delta^{-1/2}$ are Borel functions and bounded on compact sets, we always have $\langle f, p \rangle \in \mathcal{M}_{\theta c}(R)$. We must show that this mapping is determined by the equivalence class of $p$. If $p = p'$ a.e. relative to $\lambda^\mathcal{Q}$, then for $\alpha^\mu$-almost every pair $(x, y)$ the functions $p$ and $p'$ agree a.e. with respect to $\beta_y^x$, so for every $f \in \mathcal{M}_c(G)$ we have $\langle f, p \rangle = \langle f, p' \rangle$ a.e. with respect to $\alpha^\mu$. Furthermore, we represent $p$ and $p'$ as matrix entries, and these have restrictions to $X$ that agree a.e. with respect to $\mathcal{Q}$. We will show that the mapping of $f$ to $\langle f, p \rangle$ gives rise to a completely positive map $S_p$ from $M^*(G)$ to $M^*(R)$.

There is another property of $S_p$ we use, and its statement requires a little background. Recall from Section 1 that $C_c(G)$ and $\mathcal{M}_c(G)$ are bimodules over $C(\bar{X})$, where $h \in C(\bar{X})$ acts via multiplication by $h \circ r$ and $h \circ s$. Recall also that every representation $\pi$ of the $*$-algebra $C_c(G)$ has an associated representation $\varphi$ of $C_c(X)$ such that $\pi(hf) = \varphi(h)\pi(f)$ and $\pi(fh) = \pi(f)\varphi(h)$ for all $f$ and $h$, i.e., so that $\pi$ is a bimodule map. Hence every representation of $\mathcal{M}_c(G)$ also has such an associated representation of $C_c(X)$. We can extend $\varphi$ to $\mathcal{M}(X)$, getting a representation that preserves monotone limits and hence maintaining the bimodule property.

We notice that $\mathcal{M}(X)$ also has natural actions defined the same way on $\mathcal{M}_{\theta c}(R)$ and by pointwise multiplication on each $L^2(\mu; \mathcal{K})$, rendering $\tilde{\omega}$ a bimodule map from $\mathcal{M}_{\theta c}(R)$ to $M^*(R)$. The main properties of $S_p$ are established in the next theorem.

**Theorem 4.5.** If $p \in \mathcal{P}(G)$, there is a completely positive operator $S_p : M^*(G) \to M^*(R)$ that extends the operator defined by $S_p(\omega(f)) = \tilde{\omega}(\langle f, p \rangle)$ for $f \in \mathcal{M}_c(G)$. This mapping is an $\mathcal{M}(X)$-bimodule map. If we define $S_p(\omega(g\epsilon)) = \tilde{\omega}(pg\epsilon)$ for $g \in \mathcal{M}(X)$ and use linearity, we get an extension of the original $S_p$ to a completely positive $\mathcal{M}(X)$-bimodule map of $M^*(G, \bar{X})$ to $M^*(R, \bar{X})$ that takes $\omega(\epsilon)$ to an element of $\tilde{\omega}(\mathcal{M}(X))$. The completely bounded norm of $S_p$ is equal to $\|p\|_\infty$.

**Proof.** We need another formula for $S_p$, first on $\mathcal{M}_c(G, \bar{X})$. To find one, we first work with a subrepresentation of $\tilde{\omega}$ acting on a space of the form $L^2(\mu; \mathcal{K})$.

The positive definite function $p$ determines a unitary representation $\pi_p$ of $G$ on a Hilbert bundle $\mathcal{K}_p$ over $X$, as well as a bounded section $\xi_p$ of $\mathcal{K}_p$ for which we have $p(\gamma) = (\pi_p(\gamma)\xi_p \circ s(\gamma) \mid \xi_p \circ r(\gamma))$ for almost all $\gamma$ relative to $\lambda^\mathcal{Q}$. Then we may replace $p$ by the matrix entry. Indeed, we must make that replacement in order to make sense of the values of $p$ on $X$. Suppose that $\xi$ and $\eta$ are in $L^2(\mu, \mathcal{K})$, and compute:

$$(S_p(\omega(f))\xi \mid \eta) = (\tilde{\omega}(\langle f, p \rangle)\xi \mid \eta)$$

$$= \int \langle f, p \rangle(x, y)(\tilde{\omega}(x, y)\xi(y) \mid \eta(x)) \, d\tilde{\nu}_0(x, y)$$

$$= \int \int f(\gamma)p(\gamma)\delta(\gamma)^{-1/2}(\tilde{\omega} \circ \theta(\gamma)\xi(y) \mid \eta(x)) \, d\beta_y^x(\gamma) \, d\tilde{\nu}_0(x, y)$$

$$= \int f(\gamma)(\pi_p(\gamma)\xi_p \circ s(\gamma) \mid \xi_p \circ r(\gamma))(\tilde{\omega} \circ \theta(\gamma)\xi \circ s(\gamma) \mid \eta \circ r(\gamma)) \, d\nu_0(\gamma)$$

$$= (((\pi_p \otimes \tilde{\omega} \circ \theta)(f))\xi_p \otimes \xi \mid \xi_p \otimes \eta).$$

We also have

$$(\tilde{\omega}(pg\epsilon)\xi \mid \eta) = (pg\xi \mid \eta)$$

$$= \int (\xi_p(x) \mid \xi_p(x))g(x)(\xi(x) \mid \eta(x)) \, d\mu(x)$$

$$= ((\pi_p \otimes \tilde{\omega} \circ \theta)(g\epsilon)\xi_p \otimes \xi \mid \xi_p \otimes \eta).$$

Now define $V_{p,\mu,\mathcal{K}} : L^2(\mu; \mathcal{K}) \to L^2(\mu; \mathcal{K}_p \otimes \mathcal{K})$ by $V_{p,\mu,\mathcal{K}}\xi = \xi_p \otimes \xi$ and let $V$ be the direct sum of all the operators $V_{p,\mu,\mathcal{K}}$. The calculations just done show that for all $f \in \mathcal{M}_c(G)$ and $g \in \mathcal{M}(X)$ we have

$$S_p(\omega(f\lambda + g\epsilon)) = V^*((\pi_p \otimes \tilde{\omega} \circ \theta)(f\lambda + g\epsilon))V.$$

Since $\pi_p \otimes \tilde{\omega} \circ \theta$ is a $*$-representation, Stinespring's Theorem [P,St] shows that $S_p$ is completely positive. This representation also gives a formula for the extension of $S_p$ to $M^*(G,X)$ and shows that it is an extension by continuity. It is not difficult to show that the norm of $S_p$ is the essential supremum norm of $\xi_p$, and that is the same as $\|p\|_\infty$.

From the definition of $V$ we see that it intertwines the natural actions of $\mathcal{M}(X)$ on $L^2(\mu; \mathcal{K})$ and $L^2(\mu; \mathcal{K}_p \otimes \mathcal{K})$. The restrictions of these natural actions to $C_c(X)$ are the representations of $C_c(X)$ associated with the given representations of $C_c(G)$ in the proof of Renault's Theorem. This makes it clear that $S_p$ is also a bimodule map.

Now we want to prepare the way for the proof of the converse of the last theorem. We need less hypothesis than we had conclusion, namely we only need to deal with the transitive quasiinvariant measures on $X$.

We use the measures $\mu^x$ on $X$ such that $\alpha^x = \epsilon^x \times \mu^x$, as described in Section 2. For each $x$ we have $\alpha^{\mu^x} = \mu^x \times \mu^x$, which is symmetric, so $\tilde{\Delta}_{\mu^x}$ is trivial. That means that $\Delta_{\mu^x} = \delta$. Since these modular functions are all the same, we will denote them by the single letter $\Delta$.

Let $\rho_x$ be the representation of $I(R,\alpha)$ gotten by integrating the trivial representation of $R$ on the one-dimensional bundle, relative to the measure $\mu^x$. Since $\mathcal{M}_{\theta c}(R) \subseteq I(R,\alpha)$, the representation $\rho_x$ can be restricted to $\mathcal{M}_{\theta c}(R)$, and we denote the restriction the same way. Define $\rho_x$ on $\mathcal{M}(X)$ to be the representation by multiplication on $L^2(\mu^x)$. We combine these two definitions to get a representation $\rho_x$ of $\mathcal{M}_{\theta c}(R, \bar{X})$ on $\mathcal{H}_x$. Let $\tilde{\omega}_t$ denote the direct sum of all these 'transitive' representations $\rho_x$, so the representation space of $\tilde{\omega}_t$ is $\mathcal{H}_X$, the direct sum of all the Hilbert spaces $\mathcal{H}_x$. Write $M_t^*(R, \bar{X})$ for the norm closure of the image of $\mathcal{M}_{\theta c}(R, \bar{X})$ under $\tilde{\omega}_t$. Then $M_t^*(R, X)$ is a quotient of $M^*(R, X)$ as a $C(\bar{X})$-bimodule, as well as a compression of $M^*(R, X)$. We also write $M_t^*(R)$ for the closure of the image of $\mathcal{M}_{\theta c}(R)$.

It is not true that every completely bounded map is a linear combination of completely positive maps, unless the range algebra is injective [P; Notes to Chapter 7]. In our setting, the domain and range are closely related and very special. We can circumvent the problems caused by lack of injectivity, but to do so and even to deal with completely positive maps themselves, we need to think of $M_t^*(R, X)$ as acting on a space of Borel sections. We now begin to arrange that.

Observe that the Hilbert spaces $\mathcal{H}_x$ are the fibers in a Hilbert bundle over $X$, i.e., the graph of $\mathcal{H}$, $\Gamma_\mathcal{H}$, has a natural Borel structure with all the necessary properties. In fact, for each $x$ the space $\mathcal{H}_x$ is easily identifiable with $L^2(\alpha^x)$, and we simply transport the usual Borel structure for the latter bundle to $\mathcal{H}$.

If $g \in \mathcal{M}_{\theta c}(R)$, define a section of $\Gamma_\mathcal{H}$ by letting $\xi_g(x)$ be the class of $g(x, \cdot)$ in $L^2(\mu^x)$. Countably many of these sections can be chosen so that their values at a point $x$ always form a dense set in $\mathcal{H}_x$. Thus we can also choose a countably generated subalgebra of $\mathcal{M}(X)$ so that the module of sections over it generated by the countably many $\xi_g$'s determines the Borel structure on $\Gamma_\mathcal{H}$. Note also that $x \sim y$ implies that $\mu^x = \mu^y$ so $\mathcal{H}_x = \mathcal{H}_y$.

**Theorem 4.6.** Let $\psi$ be a completely positive $C(\bar{X})$-bimodule map from $C^*(G, \bar{X})$ to $M^*(R, X)$, and suppose that $\psi(\omega(\epsilon))|\mathcal{H}_X \in \tilde{\omega}_t(\mathcal{M}(X))$. Then there is a $p \in \mathcal{P}(G)$ such that $\psi = S_p$, and $\|p\|_\infty \leq \|\psi\|_{c.b.}$.

**Proof.** There is no loss of generality in taking $\psi$ to have completely bounded norm at most 1. Next we restrict $\psi \circ \omega$ to $C_c(G, \bar{X})$, getting a completely positive map, $\psi'$, of $C_c(G, \bar{X})$ into $M^*(R, X)$. For each $x \in X$, $f \in C_c(G)$, and $g \in C(\bar{X})$, define $\psi'_x(f\lambda + g\epsilon) = \psi'(f\lambda + g\epsilon)|\mathcal{H}_x$. For each $x$, $\psi'_x$ is a completely positive bimodule map into $\mathcal{L}(\mathcal{H}_x)$ of completely bounded norm at most 1, and $x \sim y$ implies $\psi'_x = \psi'_y$.

**Outline of the proof:**

The proof consists mainly of accumulating sufficient information about the mappings $\psi'_x$ and objects constructed from them to assemble the desired positive definite function $p$. Using the Stinespring Theorem for completely positive maps and analyzing the equipment it provides enables us to show that each $\psi'_x$ is of the form $S_{p_x}$. Then it is necessary to merge the separate $p_x$'s into one $p$, using the fact that $x \sim x'$ implies $\psi'_x = \psi'_{x'}$ from which we prove that $p_x = p_{x'}$ a.e. Several more improvements in the behavior of the functions $p_x$ finally allow us to produce a matrix entry that serves as the desired function $p$. We hope that naming the major steps in the proof will help the reader maintain some sense of the organization of the proof.

**Step 1: The Borel Behavior of $x \mapsto \psi'_x$.**

If $f, h \in \mathcal{M}_{\theta c}(R)$ we want to see that
$$x \mapsto \rho_x(f)(\xi_h(x))$$
is a Borel section of $\Gamma_{\mathcal{H}}$. To do this it is sufficient to show that if $f, h, k \in \mathcal{M}_{\theta c}(R)$ then the function $x \mapsto (\rho_x(f)\xi_h(x) \mid \xi_k(x))$ is Borel. Such an inner product is given by an integral, according to the definition of $\rho_x$, namely
$$\int \int f(y,z) h(x,z) \bar{k}(x,y) \, d\mu^x(z) d\mu^x(y).$$
This integral defines a Borel function of $x$ since the measures $\mu^x \times \mu^x$ depend on $x$ in a Borel manner. By the definition of $M_t^*(R)$, every $\rho_x$ is defined on $M_t^*(R)$ and for $a \in M_t^*(R)$ the function $x \mapsto \rho_x(a)$ is a uniform limit of functions of the form $x \mapsto \rho_x(f)$ for $f \in \mathcal{M}_{\theta c}(R)$. Hence for $a \in M_t^*(R)$ and $h \in \mathcal{M}_{\theta c}(R)$ the section $x \mapsto \rho_x(a)(\xi_h(x))$ is Borel.

If we define $\hat{\psi}$ to be the direct sum of all the $\psi'_x$'s, then $\hat{\psi}$ is also the compression of $\psi'$ to $\mathcal{H}_X$. Thus $\hat{\psi}$ maps $C_c(G, \bar{X})$ into $M_t^*(R, X)$ and $\rho_x \circ \hat{\psi} = \psi'_x$. From this it follows that if $f \in C_c(G)$ and $h \in \mathcal{M}_{\theta c}(R)$ then the section $x \mapsto \psi'_x(f)(\xi_h(x))$ of $\Gamma_{\mathcal{H}}$ is Borel. If $g \in C(\bar{X})$ there is a function $g_1 \in \mathcal{M}(X)$ such that $\hat{\psi}(g\epsilon) = \tilde{\omega}_t(g_1)$ because $\hat{\psi}(\epsilon) \in \tilde{\omega}_t(\mathcal{M}(X))$ and $\hat{\psi}$ is a $C(\bar{X})$-bimodule map. Hence for $a \in C_c(G, \bar{X})$ and $h \in \mathcal{M}_{\theta c}(R)$ the section $x \mapsto \psi'_x(a)(\xi_h(x))$ is Borel.

The fact that $\hat{\psi}$ maps into $M_t^*(R, X)$, and the Borel property derived above are essential for completing the proof.

**Step 2: The Stinespring Construction.**

For each $x$ we represent $\psi'_x$ by Stinespring's Theorem: We get a representation $\pi_x$ of $C_c(G, \bar{X})$ on a Hilbert space $\mathcal{K}_x$ and an operator $V_x$ from $\mathcal{H}_x$ to $\mathcal{K}_x$, such that for $a \in C_c(G, \bar{X})$ we have
$$\psi'_x(a) = V_x^* \pi_x(a) V_x.$$

We will use the details of the construction, so we repeat it here. For Stinespring's proof, it suffices to have the domain of the completely positive map to be a $*$-algebra with identity, so $C_c(G, \bar{X})$ can be used. The space $\mathcal{K}_x$ is taken to be the Hilbert space constructed from the algebraic tensor product $C_c(G, \bar{X}) \otimes \mathcal{H}_x$ using the semi-inner product whose value on two elementary tensors is given by $(a \otimes \xi \mid b \otimes \eta) = (\psi'_x(b^*a)\xi \mid \eta)$. Let $q_x$ be the quotient map from $C_c(G, \bar{X}) \otimes \mathcal{H}_x$ to its quotient modulo vectors of norm 0. The image of $q_x$ is identified with a dense subspace of $\mathcal{K}_x$. (Since $C_c(G, \bar{X})$ and $\mathcal{H}_x$ are separable, so is $\mathcal{K}_x$.) The representation $\pi_x$ is determined by having $\pi_x(a)(q(b \otimes \xi)) = q_x(ab \otimes \xi)$ for $a, b \in C_c(G, \bar{X})$ and $\xi \in \mathcal{H}_x$. The operator $V_x$ is determined by setting $V_x(\xi) = q_x(1 \otimes \xi)$ for $\xi \in \mathcal{H}_x$. A calculation of inner products shows that $\|V_x\|^2 = \|\psi'_x(1)\|$.

Since $\psi_x$, $\pi_x$, $\mathcal{K}_x$, and $V_x$ are Borel in $x$ and constant on equivalence classes, we get a Hilbert bundle over $X$ that is constant on equivalence classes. The pair $(\pi_x, V_x)$ is minimal in the sense that $\pi_x(C_c(G, \bar{X}))V_x(\mathcal{H}_x)$ is dense in $\mathcal{K}_x$.

**Step 3: Getting $p_x$ from the Stinespring Representation.**

Now we study this structure for a fixed $x \in X$. By Theorem 1.2 we know that $\pi_x$ can be obtained by integrating a representation, $\pi'_x$, of $G$ on a bundle $\mathcal{K}^x$ relative to a quasiinvariant measure $\mu_x$, i.e., $\mathcal{K}_x = L^2(\mu_x; \mathcal{K}^x)$. Let $\varphi_x$ be the representation of $C(\bar{X})$ on $\mathcal{K}_x$ associated with $\pi_x$ as a representation of $C_c(G, \bar{X})$. In terms of the representation of $\mathcal{K}_x$, $\varphi_x$ is the natural representation by multiplication on sections of $\mathcal{K}^x$ [Re2]. We also have $\varphi_x = \pi_x|C(\bar{X})$, where $C(\bar{X})$ is regarded as a subalgebra of $C_c(G, \bar{X})$. We denote the natural representation of $C(\bar{X})$ on $\mathcal{H}_x$ by $\theta_x$; again this is a representation by multiplication.

We need to show that $\mu_x$ can be taken to be $\mu^x$. The first step is to show that $V_x$ intertwines $\theta_x$ and $\varphi_x$. Take $h \in C(\bar{X})$, $b \in C_c(G, \bar{X})$, and $\xi, \eta \in \mathcal{H}_x$. Then the definition of the inner product and the fact that $\psi'_x$ is a $C(\bar{X})$-bimodule map gives

$$\begin{aligned}(1 \otimes h\xi \mid b \otimes \eta) &= (\psi'_x(b^*)h\xi \mid \eta) \\ &= (\psi'_x(b^*h)\xi \mid \eta). \\ &= (h \otimes \xi \mid b \otimes \eta)\end{aligned}$$

Hence $q_x(h \otimes \xi) = q_x(1 \otimes h\xi)$. Using the bimodule property of $\psi'_x$, the definition of $\pi_x$, and the inner product on $\mathcal{K}_x$, we compute that

$$\begin{aligned}(V_x(h\xi) \mid q_x(b \otimes \eta)) &= (q_x(h \otimes \xi) \mid q_x(b \otimes \eta)) \\ &= (\pi_x(h\epsilon)q_x(1 \otimes \xi) \mid q_x(b \otimes \eta)).\end{aligned}$$

Hence, $V_x(h\xi) = \varphi_x(h)V_x(\xi)$.

From the theory of representations of $C_c(X)$ or of projection valued measures based on $X$, there is a bounded section of $\mathcal{K}^x$, which we denote by $\zeta_x$, such that for $\xi \in \mathcal{H}_x$, the pointwise product $\xi\zeta_x$ is a section of $\mathcal{K}^x$ representing the element $V_x(\xi)$ in $\mathcal{K}_x$. Such a section can be gotten as follows: let $g$ be any strictly positive Borel function on $X$ that represents an element of $\mathcal{H}_x$, let $\zeta^1$ be a section that represents $V_x(g)$, and set $\zeta_x = (1/g)\zeta^1$. Then $\zeta_x$ need not be a square integrable section, but will be if $\mu^x$ is finite so that the function 1 is an element of $\mathcal{H}_x$.

We can write $V_x(\xi) = \xi\zeta_x$, using the usual identification of functions with their equivalence classes. Then for $\xi \in \mathcal{H}_x$ we have

$$(*) \qquad \int |\xi|^2 |\zeta_x|^2 \, d\mu_x \leq \int |\xi|^2 \, d\mu^x,$$

because $\|V_x\| \leq 1$. It follows that $\mu_x$ is not singular relative to $\mu^x$, so that $\mu_x$ gives positive measure to $[x]$. It also follows that $|\zeta_x|$ is zero a.e. off $[x]$, so that $\zeta^1 = g\zeta_x$ is in the subspace of $\mathcal{K}_x = L^2(\mu_x; \mathcal{K}^x)$ consisting of functions that vanish off $[x]$. By the way we integrate representations of $G$ to get representations of $C_c(G)$, we see that this latter subspace is invariant for $C_c(G)$ and hence for $C_c(G, \bar{X})$. From the fact that $g$ is cyclic in $\mathcal{H}_x$, it follows that $g\zeta_x = V_x(g)$ is cyclic for $C_c(G, \bar{X})$ in $\mathcal{K}_x$, so the subspace under discussion is in fact all of $\mathcal{K}_x$. That implies that $\mu_x$ is in fact equivalent to $\mu^x$, so we may as well take $\mu_x$ to be equal to $\mu^x$. That may require multiplying the original $\zeta_x$ by some positive function, but now we assume that to have been done. We write $\nu^x$ for $\lambda^{\mu^x}$, getting a measure concentrated on $G|[x]$.

In this situation, the inequality $(*)$ implies that $|\zeta_x|$ is bounded by 1. We define

$$p_x(\gamma) = (\pi'_x(\gamma)\zeta_x(s(\gamma)) \mid \zeta_x(r(\gamma)))$$

getting a positive definite function on $G|[x]$. Now the sup-norm of $\zeta_x$ is the same as the operator norm of $V_x$, and that is the same as the square root of the completely bounded norm of $\psi'_x$, so the sup-norm of $p_x$ is at most the completely bounded norm of $\psi'_x$.

**Step 4.** $p_x$ **gives rise to** $\psi'_x$. We know that $x \sim y$ implies $\psi'_x = \psi'_y$, so $\pi_x = \pi_y$ and $V_x = V_y$. Hence $\pi'_x(\gamma) = \pi'_y(\gamma)$ for $\nu^x$-almost every $\gamma$, and $\zeta_x(z) = \zeta_y(z)$ for $\mu^x$-almost every $z$, so that $p_x = p_y$ a.e. relative to $\nu^x$, and their restrictions to $X$ agree a.e. relative to $\mu^x$.

To see that $\psi'_x$ is the compression of $S_{p_x}$ to $\mathcal{H}_x$, we begin by setting $\nu^x = \lambda^{\mu^x}$ and $\tilde{\nu}^x = \alpha^{\mu^x}$, as above, so that $\tilde{\Delta} = 1$ and $\Delta = \delta$. Then we calculate for $f \in C_c(G)$, and $\xi, \eta \in \mathcal{H}_x$:

$$\begin{aligned}(\psi'_x(f)\xi \mid \eta) &= (\pi_x(f)V_x\xi \mid V_x\eta) \\ &= \int f(\gamma)(\pi_x(\gamma)(\xi\zeta_x)(s(\gamma)) \mid (\eta\zeta_x)(r(\gamma)))\Delta^{-1/2}(\gamma) \, d\nu^x(\gamma) \\ &= \int\int f(\gamma)p_x(\gamma)\delta^{-1/2}(\gamma) \, d\beta^y_z(\gamma) \, \xi(z)\bar{\eta}(y) \, d\tilde{\nu}^x(y,z).\end{aligned}$$

This shows that $\psi'_x(f) = S_{p_x}(f)|\mathcal{H}_x$. Next we find a formula for $\psi'_x(\epsilon)$ by computing

$$(\psi'_x(\epsilon)\xi \mid \eta) = (\pi_x(\epsilon)V_x\xi \mid V_x\eta)$$
$$= (\xi\zeta_x \mid \eta\zeta_x)$$
$$= \int p_x(y)\xi(y)\bar{\eta}(y)\,d\mu^x(y)$$

from which it follows that $\psi'_x(\epsilon) = \rho_x(p_x|X)$. Since $\psi'_x$ is a $C(\bar{X})$-bimodule map, we see that $\psi'_x(g\epsilon) = \rho_x(gp_x)$ for $g \in C(\bar{X})$. This completes the proof that $\psi'_x$ is the compression of $S_{p_x}$ to $\mathcal{H}_x$.

**Step 5. Applying Lemma 2.1 to the Functions $p_x$.**

Take functions $h, k \in \mathcal{M}_{\theta c}(R)$ from which we make sections $\xi_h$ and $\xi_k$ of $\mathcal{H}$. Let $\xi = \xi_h(x)$ and $\eta = \xi_k(x)$ in the calculations above to see that if $g \in C_c(G)$, then

$$(\psi'_x(g)\xi_h(x) \mid \xi_k(x)) = \int g(\gamma)p_x(\gamma)h(x, s(\gamma))\bar{k}(x, r(\gamma))\delta^{-1/2}(\gamma)\,d\nu^x(\gamma).$$

If $\epsilon$ is the identity in $C(\bar{X})$, we also get

$$(\psi'_x(\epsilon)\xi_h(x) \mid \xi_k(x)) = \int p_x(\gamma)h(x, y)\bar{k}(x, y)\,d\mu^x(y).$$

Here it is important that the functions of $x$ on the left hand sides of these two formulas are Borel functions.

To apply Lemma 2.1 as it is formulated, we must have a Borel family of finite measures. We begin by considering a compact set $K$ contained in $G$. The function $y \mapsto \lambda^y(K)$ is bounded on $X$, and for every $x \in X$ we have $\mu^x(s(xK)) < \infty$. Hence, for $x \in X$ the measure given by the integral

$$\int_{s(xK)} (\chi_K \lambda^y)\,d\mu^x(y)$$

is finite.

Notice that a pair $(x, y) \in X \times X$ is in $\theta(K)$ iff $x \in r(Ky)$ iff $y \in s(xK)$. If $h$ is the characteristic function of $\theta(K)$, it follows that $h(x, r(\gamma)) = 1$ iff $r(\gamma) \in s(xK)$, and $h(x, s(\gamma)) = 1$ iff $s(\gamma) \in s(xK)$. Thus the set $L$, defined to be $\{(x, \gamma) \in X \times G : \gamma \in K \text{ and } h(x, s(\gamma))h(x, r(\gamma)) = 1\}$ is a Borel set in $X \times G$, and the same as $\{(x, \gamma) \in X \times G : \gamma \in K, s(\gamma) \in s(xK) \text{ and } r(\gamma) \in s(xK)\}$. From the preceding paragraph, it follows that every $x$-section $L_x$ of $L$ has finite measure for $\nu^x$.

Choose compact sets $K_1 \subset K_2 \subset \cdots$ whose union is $G$, and for each $n$ define $h_n = \chi_{\theta(K_n)}$ and then $L_n = \{(x, \gamma) \in X \times G : \gamma \in K_n, s(\gamma) \in s(xK_n) \text{ and } r(\gamma) \in s(xK_n)\}$. Define $D_1 = L_1$ and for $n \geq 2$, let $D_n = L_n \setminus L_{n-1}$. For each $n \in \mathbb{N}$ and $x \in X$, let $\nu_n^x = (\chi_{(D_n)_x})\nu^x$. This gives a Borel family of finite measures on $G$. Notice that the sets $D_n$ partition $\{(x, \gamma) \in X \times G : \gamma \in G|[x]\}$.

Now define $f_x$ on $G$ for $x \in X$ by $f_x(\gamma) = p_x(\gamma)\delta^{-1/2}(\gamma)$ for $\gamma \in G|[x]$ and $0$ for other $\gamma$'s. If $g \in C_c(G)$ and $x \in X$, then

$$\int g(\gamma)f_x(\gamma)\,d\nu_n^x(\gamma) = (\psi'_x(g)\xi_{h_n}(x) \mid \xi_{h_n}(x)),$$

which is a Borel function of $x$. Hence there is a Borel function $F_n$ on $X \times G$ such that for each $x$, $F_n(x, \cdot) = f_x$ a.e. relative to $\nu_n^x$. Set

$$F = \sum_{n \geq 1} \chi_{D_n} F_n.$$

Then $F$ is Borel and for each $x \in X$, $F(x, \cdot) = f_x$ a.e. relative to $\nu^x$.

A similar analysis using $\mu^x$ shows that we can also choose $F$ so that $F(x, y) = f_x(y)$ for $\mu^x$-almost every $y$.

Hence there is a Borel function $P$ on $X \times G$ such that for every $x$ we have $P(x, \cdot) = p_x$ a.e. Also, $x \sim y$ implies that $P(x, \cdot) = P(y, \cdot)$ a.e. relative to either $\nu^x$ or $\nu^y$ (these are the same measure) and also relative to either $\mu^x$ or $\mu^y$ when restricted to $X$. Furthermore, $|P(x, \cdot)|$ is bounded by the completely bounded norm of $\psi'_x$, so $|P|$ is bounded by 1.

**Step 6. Improving the behavior of $P$.**

Recall the probability measures $\mu_1^x = s(\lambda_1^x)$ on $X$ obtained from the Borel family of normalized Haar measures on $G$. We have $\mu_1^x \sim \mu_1^y$ if $x \sim y$. Define a new function $P_1$ on $X \times G$ by

$$P_1(x, \gamma) = \int P(y, \gamma) \, d\mu_1^x(y).$$

Make a function of three variables from $P$ and use the Borel character of $P$ and the measures $\mu_1^x$ to show that $P_1$ is also Borel. We need to know that $P_1$ also essentially replicates every function $p_x$, and is even more invariant than $P$ under changing $x$ to an equivalent point of $X$.

To begin with we limit ourselves to one orbit, and denote it by $S$. We write $\mu^S$ for a choice of one of the measures $\mu_1^{x_0}$ for $x_0 \in S$. We know that for $x$ and $y$ in $S$ the functions $P(x, \cdot)$ and $P(y, \cdot)$ agree a.e. relative to $\lambda^{\mu^S}$, so they agree a.e. relative to $\lambda^z$ for $\mu^S$-almost every $z$. Since $\lambda^z$ and $\lambda_1^z$ have the same null sets, $P(x, \cdot)$ and $P(y, \cdot)$ agree a.e. relative to $\lambda^z$ iff the complex measures $P(x, \cdot)\lambda_1^z$ and $P(y, \cdot)\lambda_1^z$ are the same. We have two Borel mappings from $S^3$ to the standard Borel space of complex Borel measures on $\bar{X}$, so the set $E_S$ on which they agree is Borel, allowing us to use Fubini arguments.

Hence, for every $x \in S$, the set $\{(y, z) \in S^2 : P(y, \cdot) = P(x, \cdot)$ a.e. $d\lambda^z\}$ is a Borel set whose complement has measure 0 for $\mu^S \times \mu^S$. Therefore, there is a conull Borel set $Z_x$ of points $z$ in $S$ such that for $\mu^S$-almost every $y$ we have $P(y, \cdot) = P(x, \cdot)$ a.e. relative to $\lambda^z$. Thus, for $z \in Z_x$ it is true that for $\lambda^z$-almost every $\gamma$ we have $P(y, \gamma) = P(x, \gamma)$ for $\mu^S$-almost every $y$. It follows that if $z \in Z$, then $P_1(x, \gamma) = P(x, \gamma)$ for $\lambda^z$-almost every $\gamma$. Hence, for every $x \in S$ we have $P_1(x, \cdot) = P(x, \cdot)$ a.e. In particular, $P_1$ also replicates every $p_x$, since $S$ is a general orbit.

In the last paragraph, we enountered points $\gamma \in G$ for which $P(y, \gamma)$ is essentially constant in $y$ because it is almost always equal to a particular $P(x, \gamma)$. We need to know more about the set $H = \{\gamma \in G : y \mapsto P(y, \gamma)$ is essentially constant$\}$. If $\mathcal{A}$ is a countable algebra that generates the Borel sets in $\mathbb{C}$, it is not difficult to show that

$$H = \bigcap_{A \in \mathcal{A}} \{\gamma \in G : \mu_1^{r(\gamma)} \times \epsilon^\gamma(P^{-1}(A)) \in \{0, 1\}\}.$$

Thus $H$ is a Borel set. Hence the set $C = \{x \in X : \lambda_1^x(H) = 1\}$ is also a Borel set. From the preceding paragraph, it follows that $C$ is conull in every orbit. For $z \in C$, the function $P_1(\cdot, \gamma)$ is constant for $\lambda^z$-almost every $\gamma \in zG$. In particular, for $z \in C$ it is true that $x, y \in [z]$ implies that $P_1(x, \cdot)\lambda_1^z = P(y, \cdot)\lambda_1^z$.

The last conclusion is the additional invariance needed, and now we change notation and simply write $P$ for $P_1$, since it does everything we need.

**Step 7. Making a Borel family of representations from $P$.**

Again, take a particular orbit, $S$, in $X$. For every pair $(x, y) \in S^2$, we have $P(x, \cdot) = P(y, \cdot)$ a.e. relative to $\lambda^z$ for $\mu^S$-almost every $z$. Take an arbitrary $z \in S$. Then for $\lambda^z$-almost every $\gamma_2$ it is true that $P(x, \cdot) = P(y, \cdot)$ a.e. relative to $\gamma_2^{-1} \cdot \lambda^z = \lambda^{s(\gamma_2)}$. Hence $P(x, \gamma_2^{-1}\gamma_1) = P(y, \gamma_2^{-1}\gamma_1)$ for $\lambda^z \times \lambda^z$-almost every pair $(\gamma_1, \gamma_2)$. (The mapping taking the pair to $\gamma_2^{-1}\gamma_1$ carries $\lambda_1^z \times \lambda_1^z$ to a measure equivalent to $\lambda^{\mu^z}$.)

Now return to studying general points of $X$. For $f, g \in C_c(G)$ and $(x, y) \in R$, define

$$(f \mid g)_{(x,y)} = \int \int f(\gamma_1) \bar{g}(\gamma_2) P(x, \gamma_2^{-1}\gamma_1) \, d\lambda^y(\gamma_1) d\lambda^y(\gamma_2).$$

The formula defines an inner product on $C_c(G)$, and we write $\mathcal{K}(x, y)$ for the resulting Hilbert space. For each $f, g \in C_c(G)$ the function $(x, y) \mapsto (f, g)_{(x,y)}$ is a Borel function on $R$ that is constant on sets of the form $[y] \times \{y\}$, so $\mathcal{K}$ defines a Hilbert bundle on $R$ that is constant on the same sets. For $f \in C_c(G)$, let $\sigma(f)$ denote the section of $\mathcal{K}$ (or $\Gamma_\mathcal{K}$) that it determines.

For each $x$, the bundle $\mathcal{K}(x, \cdot)$ supports a unitary representation: here we denote it by $\pi_x$ rather than $\pi_{P(x,\cdot)}$ (see Section 3). We know that $x \sim x'$ implies that $\pi_x = \pi_{x'}$, which means that for $\gamma \in G|[x]$ we have $\pi_x(\gamma) = \pi_{x'}(\gamma)$ (they are on the same space). We want to show that $(x, \gamma) \mapsto \pi_x(\gamma)$ is Borel on $X \times' G = \{(x, \gamma) : \gamma \in G|[x]\}$. It will help to look at $R \times' G = \{(x, y, \gamma) : \gamma \in G|[x]\}$. The function

$$(x, y, \gamma) \mapsto \int \int f(\gamma^{-1} \gamma_1) \bar{g}(\gamma_2) P(x, \gamma_2^{-1} \gamma_1) \, d\lambda^y(\gamma_1) d\lambda^y(\gamma_2)$$

is Borel on $R \times' G$, so $(x, \gamma) \mapsto (\pi_x(\gamma) \sigma(f)(x, s(\gamma)) \mid \sigma(g)(x, r(\gamma)))$ is Borel on $X \times' G$.

**Step 8. Finding a Borel section that represents $P$.**

Let $D$ be the set of pairs $(x, y) \in R$ for which the linear functional $f \mapsto \lambda^y(fP(x, \cdot))$ is bounded relative to the seminorm $\|\sigma(f)(x, y)\|$ on $C_c(G)$. The boundedness can be tested using a countable dense subset of $C_c(G)$, so $D$ is Borel, and hence so is the set $DC$. For each $x \in X$, we have $xD = \{x\} \times D_x$ so that $xD$ is conull with respect to $\alpha^x$. Notice that $w \sim x$ implies that $C \cap D_x = C \cap D_w$, and this set is conull in the orbit. Hence $xDC$ and $wDC$ have the same conull image in $[x]$ under $s$.

Now, for $(x, y) \in D$ define $\zeta(x, y)$ to be the vector in $\mathcal{K}(x, y)$ such that $(\sigma(f)(x, y) \mid \zeta(x, y)) = \lambda^y(fP(x, \cdot))$ for every $f \in C_c(G)$, and for $(x, y) \notin D$, let $\zeta(x, y) = 0$. The formula makes it clear that $\zeta$ is Borel.

If $y \in C$ and $w \sim x \sim y$, then $(w, y) \in D$ iff $(x, y) \in D$, so $y \in C$ implies that $Dy = [y] \times \{y\}$. Also, $w, x \in [y]$ implies that $P(w, \cdot)$ and $P(x, \cdot)$ agree a.e. with respect to $\lambda^y$ and that $\mathcal{K}(w, y) = \mathcal{K}(x, y)$. Together, these imply that $\zeta(w, y) = \zeta(x, y)$. Then for every $\gamma \in G|[y]$,

$$(\pi_x(\gamma)\zeta(x, s(\gamma)) \mid \zeta(x, r(\gamma)) = (\pi_w(\gamma)\zeta(w, s(\gamma)) \mid \zeta(w, r(\gamma)).$$

Thus both of these functions agree a.e. on $G|[x]$ with $P(x, \cdot)$. Thus we can define

$$p(\gamma) = (\pi_{s(\gamma)}(\gamma)\zeta(s(\gamma), s(\gamma)) \mid \zeta(s(\gamma), r(\gamma)))$$

for $\gamma \in \theta^{-1}(DC)$ and 1 for other $\gamma$'s to get a Borel function on $G$ that agrees a.e. with $P(x, \cdot)$ on $G|[x]$.

From Step 4 it follows that $\hat{\psi}$ and the compression of $S_p$ to $\mathcal{H}_X$ are the same.

**Step 9. The compression map from $\mathcal{L}(\mathcal{H}_\omega)$ to $\mathcal{L}(\mathcal{H}_X)$.**

To complete the proof, need to show that the compression map $C$ from $\mathcal{L}(\mathcal{H}_\omega)$ to $\mathcal{L}(\mathcal{H}_X)$ is one-one when restricted to $\tilde{\omega}(\mathcal{M}_{\theta c}(R, \bar{X}))$. Then it will follow that $\psi$ and $S_p$ agree on $C_c(G, \bar{X})$, forcing them to be the same.

Suppose that $f\alpha + g\epsilon \in \mathcal{M}_{\theta c}(R, X)$ and $\tilde{\omega}(f\alpha + g\epsilon) \neq 0$. Then there is a representation $\pi$ of $R$ and a probability measure $\mu \in \mathcal{Q}$ such that $\pi^\mu(f\alpha + g\epsilon) \neq 0$. We need to use this to find a $z \in X$ such that $\rho_z(f\alpha + g\epsilon) \neq 0$, which will imply $\tilde{\omega}_t(f\alpha + g\epsilon) \neq 0$. There is no loss of generality in assuming that there is a probability measure $\mu'$ on $X$ such that $\mu = \int \mu_1^x \, d\mu'(x)$. Set $A = \{(x, y) \in R : x \neq y, \text{ and } f(x, y) \neq 0\}$, and consider two cases: $\alpha^\mu(A) = 0$ and $\alpha^\mu(A) \neq 0$. In the first case, $\pi^\mu(f) = 0$ unless $\alpha^\mu(X) > 0$, in which case we have $f\alpha = f\epsilon$ relative to $\alpha^\mu$. Thus there is an $h \in \mathcal{M}(X)$ such that $0 \neq \pi^\mu(f\alpha + g\epsilon) = \pi^\mu(h\epsilon)$. Then $\mu(\{h \neq 0\}) > 0$ so there is a $z \in X$ such that $\mu^z(\{h \neq 0\}) > 0$, and it is easy to show that $\rho_z(h\epsilon) \neq 0$, i.e., $\rho_z(f\alpha + g\epsilon) \neq 0$. In the second case, there is a $z \in X$ such that $\alpha^{\mu^z}(A) > 0$, and we will show that $\rho_z(f\alpha + g\epsilon) \neq 0$. Recall that $\alpha^{\mu^z} = \mu^z \times \mu^z$. Set $R_0 = R \backslash \{(x, x) : x \in X\}$. Then sets of the form $(E \times F) \cap R_0$, where $E$ and $F$ are disjoint Borel sets in $X$, generate the Borel sets in $R_0$, so there must be such a pair for which

$$0 < \int_{E \times F} f \, d(\mu^z \times \mu^z) < \infty.$$

If we set $h_1 = \chi_F$ and $h_2 = \chi_E$ we get elements of $\mathcal{M}(X)$ which we think of as elements of $\mathcal{H}_z$, and then the displayed integral is $(\rho_z(f)h_1 \mid h_2)$. On the other hand, $(\rho_z(g\epsilon)h_1 \mid h_2) = 0$ because $gh_1\bar{h}_2 = 0$. Thus $\rho_z(f\alpha + g\epsilon) \neq 0$, as needed.

## 5. Completely Bounded Bimodule Maps

Recall that $\mathcal{B}(G)$ is defined to be the linear span of $\mathcal{P}(G)$. Because we know that $\mathcal{P}(G)$ consists of diagonal matrix entries of unitary representations we can form direct sums of representations to show that elements of $\mathcal{B}(G)$ are also matrix entries that need not be diagonal. In this section we will provide $\mathcal{B}(G)$ a normed algebra structure. One way to compute the norm of an element $b$ of $\mathcal{B}(G)$ is in terms of the positive definite functions on a larger groupoid for which $b$ can appear as an "off diagonal part." This is the groupoid version of the well known $2 \times 2$ matrix method, and has been exploited by Renault for the same purpose [Re4]. This permits using the completeness of $\mathcal{P}(G)$ for a general locally compact groupoid to prove the completeness of $\mathcal{B}(G)$.

We can also formulate $\mathcal{B}(G)$ as an algebra of completely bounded $C(\bar{X})$-bimodule maps on $M^*(G)$, and as a space of completely bounded $C(\bar{X})$-bimodule maps from $C^*(G, \bar{X})$ to $M^*(R, X)$. Since the completely positive elements in the latter set are all given by positive definite functions, and the completely positive bimodule maps form a complete set, we get one way to prove that $\mathcal{B}(G)$ is complete.

Recall that $\omega$ is the direct sum of all cyclic representations of $C^*(G)$. We can construct each cyclic representation as an integrated representation of $G$, and, as such, it can be taken as a representation of $\mathcal{M}_c(G)$, and we use the same notation. For each $a \in C^*(G)$, $\|\omega(a)\| = \|a\|$ is the same as $\sup\{\|\pi(a)\| : \pi$ is a cyclic representation of $C^*(G)\}$. Also recall, from Section 1, the norms $\| \ \|_{II,\mu}$ and $\| \ \|_{II}$ and their properties.

**Theorem 5.1**. If $b \in \mathcal{B}(G)$, the operator $T_b$, taking $\omega(f)$ to $\omega(bf)$ for $f \in \mathcal{M}_c(G)$, extends to a completely bounded map of $M^*(G)$ to itself and $\|T_b\|_{cb} \geq \|b\|_\mathcal{Q}$.

**Proof.** By Theorem 4.1, if $p \in \mathcal{P}(G)$ then $T_p$ is completely positive, so for $b \in \mathcal{B}(G)$ the operator $T_b$ is completely bounded. Set $M = \|b\|_\mathcal{Q}$ and suppose $0 < \alpha < 1$. Since $\alpha$ is arbitrary, the proof will be complete if we find an $f \in \mathcal{M}_c(G)$ such that $\omega(f) \neq 0$ and $\|T_b\omega(f)\| \geq M\alpha^2\|\omega(f)\|$. To find such an $f$ first notice that there is a $\mu \in \mathcal{Q}$ such that the $L^\infty(\lambda^\mu)$–norm of $b$ is greater than $M\alpha$, so there exist a $b_0 \in \mathbb{C}$ and a $\eta > 0$ such that the set $A = \{\gamma : |b(\gamma) - b_0| < \eta\}$ has positive measure for $\lambda^\mu$ and $|b_0| - \eta > M\alpha$. Then there is a compact set $C \subseteq A$ such that $\lambda^\mu(C) > 0$. We take $f = \chi_C$.

By the definition of $\| \ \|_{II}$, there is a $\mu' \in \mathcal{Q}$ such that $\|f\|_{II,\mu'} > \alpha\|f\|_{II}$. By the properties of $\| \ \|_{II}$ from Section 1, if $\pi$ is the one-dimensional trivial representation of $G$, we have $\|\pi^{\mu'}(f)\| > \alpha\|\omega(f)\|$. Now let $\sigma = \pi^\mu \oplus \pi^{\mu'}$. We have $\|\sigma(f)\| \geq \|\pi^{\mu'}(f)\| > \alpha\|\omega(f)\|$.

We can find $g_1$ and $g_2$ in $C_c(X)$, $\geq 0$, and $> 0$ on $r(C) \cup s(C)$. These can be regarded as sections of the bundle for $\pi$, and it is clear that $(\pi^\mu(f)g_1|g_2) > 0$ from the integral formula for the inner product. Thus $\pi^\mu(f) \neq 0$, so $\sigma(f) \neq 0$ and $\omega(f) \neq 0$.

Since $\sigma(b_0 f) = b_0 \sigma(f)$, it will suffice to show that $\|\sigma((b - b_0)f)\| \leq \eta\|\sigma(f)\|$, because then we get $(|b_0| - \eta)\|\sigma(f)\| \leq \|\sigma(bf)\|$, so $(|b_0| - \eta)\alpha\|\omega(f)\| \leq \|\sigma(bf)\| \leq \|\omega(bf)\| = \|T_b(\omega(f))\|$, giving the desired inequality. Now $f$ is a characteristic function, so $(b - b_0)f = ((b - b_0)f)f$. Also, $\|(b - b_0)f\|_\infty \leq \eta$, so the inequality we wanted on $\sigma$ can be obtained by applying the second inequality before Lemma 1.3 to both $\mu$ and $\mu'$. Thus the proof is complete.

Again we use the algebra $C(G, \bar{X})$ to study $\mathcal{B}(G)$, and need the one-one correspondence between its representations and those of $C_c(G)$ and hence those of $G$. We still use $\omega$ for the direct sum of all cyclic representations of $C(G, \bar{X})$, each of them given as an integrated representation of $G$. We use $\tilde{\omega}$ for the direct sum of all the cyclic representations of $\mathcal{M}_c(R, X)$ that can be obtained by integrating a representation of $R$. Recall that $C^*(G, \bar{X})$ is the operator-norm closure of $\omega(C(G, \bar{X}))$ and $M^*(R, \bar{X})$ is the operator norm closure of $\tilde{\omega}(\mathcal{M}_c(R, X))$. If $x \in X$, use $\mathcal{H}_x$ for $L^2(\mu^x)$ as before, and $\mathcal{H}_X$ for the direct sum of all the $\mathcal{H}_x$'s. Let $\tilde{\omega}_t$ be the subrepresentation of $\tilde{\omega}$ obtained by restricting to $\mathcal{H}_X$.

**Theorem 5.2.** Let $b \in \mathcal{B}(G)$. There is a completely bounded $C(\bar{X})$-bimodule map $S_b : C^*(G, \bar{X}) \to M^*(R, X)$ such that $S_b(\omega(f)) = \tilde{\omega}(<f, b>)$ for $f \in C_c(G)$ and $S_b(\omega(g\epsilon)) = \tilde{\omega}(bg\epsilon)$ for $g \in C(\bar{X})$. For this operator we have

$$\|S_b\|_{cb} \geq \|b\|_\mathcal{Q},$$

and

$$S_b(\omega(\epsilon))|\mathcal{H}_X \in \tilde{\omega}_t(\mathcal{M}(X)).$$

**Proof.** The operator $S_b$ is a linear combination of four operators $S_p$ for $p \in \mathcal{P}(G)$, and these are completely positive bimodule maps by Theorem 4.5.

For the norm inequality, we proceed as in the proof of Theorem 5.1. Let $M = \|b\|_Q$ and $0 < \alpha < 1$. It will suffice to find $f \in \mathcal{M}_c(G)$ such that $\omega(f) \neq 0$ and $\|S_b\omega(f)\| \geq M\alpha^2\|\omega(f)\|$. Choose $\mu, b_0, \eta, A$ and $C$ as in Theorem 5.1, and take $f = \chi_C$.

We take $\pi$ to be the trivial one-dimensional representation, and choose $\mu'$ and $\sigma$ as before. The proof that $\omega(f) \neq 0$ used before works here also.

Let $\tilde{\pi}$ denote the one-dimensional trivial representation of $R$, and form its integral with respect to $\mu$, $\tilde{\pi}^\mu$. Likewise form $\tilde{\pi}^{\mu'}$, and let $\tilde{\sigma} = \tilde{\pi}^\mu \oplus \tilde{\pi}^{\mu'}$. It will suffice to prove that $\|\tilde{\sigma}(\langle f, 1\rangle)\| > \alpha\|\omega(f)\|$.

For this purpose, we need to see that $\|\langle f, 1\rangle\|_{II,\mu} = \|f\|_{II,\mu}$. This follows from the fact that $f \geq 0$ together with the relationship between $\nu_0$ and $\tilde{\nu}_0$. Then we see that

$$\|\langle f, b - b_0\rangle\|_{II,\mu} \leq \|(b - b_0)f\|_Q \|f\|_{II,\mu} < \eta \|f\|_{II,\mu}$$

using the fact that $f$ is a characteristic function.

Both the equality and the inequalities also hold for $\mu'$, and since $\pi$ and $\tilde{\pi}$ are the one-dimensional trivial representations, they transfer to the corresponding equality and inequalities for $\sigma$ and $\tilde{\sigma}$. Hence

$$\begin{aligned}\|\tilde{\omega}(\langle f, b\rangle)\| &\geq \|\tilde{\sigma}(\langle f, b\rangle)\| \\ &\geq \|\tilde{\sigma}(\langle f, b_0\rangle)\| - \|\tilde{\sigma}(\langle f, b - b_0\rangle)\| \\ &\geq |b_0|\|\tilde{\sigma}(\langle f, 1\rangle)\| - \eta\|\tilde{\sigma}(\langle f, 1\rangle)\| \\ &\geq M\alpha\|\tilde{\sigma}(\langle f, 1\rangle)\| \\ &\geq M\alpha^2\|\omega(f)\|.\end{aligned}$$

In order to provide the norm on $\mathcal{B}(G)$ in a way that will be convenient for proving completeness, we introduce a way to enlarge the groupoid $G$ as it was done in [Re4]. Write $T_2$ for the transitive equivalence relation on the two element set $\{1, 2\}$, so that $T_2$ has four elements. It will be convenient to have a shorter notation for matrix coefficients: If $\pi$ is a unitary representation of $G$ and $\xi$ and $\eta$ are bounded Borel sections of the bundle $\mathcal{H}$ on which $\pi$ acts, we can write $[\pi, \xi, \eta]$ for the matrix coefficient, namely

$$[\pi, \xi, \eta](\gamma) = (\pi(\gamma)\xi \circ s(\gamma) \mid \eta \circ r(\gamma)).$$

**Theorem 5.3.** A bounded Borel function $b$ on $G$ is in $\mathcal{B}(G)$ if and only if there is a function $p' \in \mathcal{P}(G \times T_2)$ such that for $\gamma \in G$ we have $b(\gamma) = p'(\gamma, (1, 2))$. The function $b$ can be expressed as a matrix coefficient using sections of sup norm at most 1 if and only if there is an associated $p'$ that can be expressed as a diagonal matrix coefficient using a section of sup norm at most 1.

**Proof.** The proof of the first assertion will be given in terms of matrix coefficients and will include the proofs of the facts about sup norms. Let $X' = X \times \{1, 2\}$ be the unit space of $G' = G \times T_2$.

Suppose that $\pi$ is a unitary representation of $G$ on a bundle $\mathcal{H}$ and that $\xi$ and $\eta$ are Borel sections of $\mathcal{H}$ of sup norm at most 1 such that $b = [\pi, \xi, \eta]$. Define a Hilbert bundle $\mathcal{H}'$ over $X'$ by setting $\mathcal{H}'(x, i) = \mathcal{H}(x)$ for $i = 1, 2$. For $\gamma' = (\gamma, (i, j))$ in $G'$ notice that $s(\gamma') = (s(\gamma), j)$ and $r(\gamma') = (r(\gamma), i)$. That means that we can define a representation $\pi'$ of $G'$ on $\mathcal{H}'$ by $\pi'(\gamma') = \pi(\gamma)$. Define a section $\zeta'$ of $\mathcal{H}'$ by setting $\zeta'(x, i) = \eta(x)$ when $i = 1$ and $\zeta'(x, i) = \xi(x)$ when $i = 2$. Then the sup norm of $\zeta'$ is at most 1 and for every $\gamma \in G$ we have $b(\gamma) = [\pi', \zeta', \zeta'](\gamma, (1, 2))$ as required.

For the converse, suppose we begin with $\mathcal{H}'$, $\pi'$, and $\zeta'$. Then for $x \in X$ define $\mathcal{H}(x) = \mathcal{H}'(x, 1) \oplus \mathcal{H}'(x, 2)$ and set $\eta(x) = (\zeta'(x, 1), 0)$ and $\xi(x) = (0, \zeta'(x, 2))$. For $\gamma \in G$ define $\pi(\gamma)$ to take $(\xi_1, \xi_2)$ to

$$(\pi'(\gamma, (1, 1))\xi_1 + \pi'(\gamma, (1, 2))\xi_2, \pi'(\gamma, (2, 1))\xi_1 + \pi'(\gamma, (2, 2))\xi_2),$$

thus acting as a matrix by left multiplication on column vectors. The sections $\xi$ and $\eta$ have sup norm at most 1, and we have $b = [\pi, \xi, \eta]$.

## 6. $\mathcal{B}(G)$ Is a Banach Algebra

Because of the results of Sections 3, 4 and 5 we can now complete the task we set ourselves at the beginning of the paper, as indicated by the section heading. Recall that for $b \in \mathcal{B}(G)$, $T_b$ is the operator on $M^*(G)$ detemined by multiplication by $b$ on $\mathcal{M}_c(G)$, and that we sometimes work with $\mathcal{B}(G)$ as an algebra of functions, even though the elements are actually equivalence classes.

**Theorem 6.1.** $\mathcal{B}(G)$ is a Banach algebra with pointwise operations for the algebraic structure and with the norm defined by
$$\|b\| = \|T_b\|_{cb}$$
for $b \in \mathcal{B}(G)$.

**Proof.** Theorem 3.8 shows that $\mathcal{B}(G)$ is an algebra under pointwise operations, and equals $\mathcal{P}(G) - \mathcal{P}(G) + i\mathcal{P}(G) - i\mathcal{P}(G)$. Any function that is 0 for $\lambda^{\mathcal{Q}}$-almost every point of $G$ represents the 0 element of $M^*(G)$, so for $b \in \mathcal{B}(G)$ the operator $T_b$ depends only on the equivalence class of $b$. Thus $b \mapsto T_b$ is well defined from the space of equivalence classes of functions in $\mathcal{B}(G)$ to the space of completely bounded operators on $M^*(G)$. Since $\|T_b\|_{cb} \geq \|b\|_\infty$, we see that $b \mapsto T_b$ is also one-one. Thus the norm makes $\mathcal{B}(G)$ a commutative normed algebra.

To prove that $\mathcal{B}(G)$ is complete, let $b_1, b_2, \ldots$ be a sequence in $\mathcal{B}(G)$ such that the norms $\|T_{b_n}\|_{cb}$ are summable. Then Theorem 5.3 says that we can construct positive definite functions $p'_1, p'_2, \ldots$ on the groupoid $G' = G \times T_2$ of Section 5 such that for every $\gamma \in G$ and every $n$ we have $b_n(\gamma) = p'_n(\gamma, (1,2))$, and for every $n$ we have $\|p'_n\|_\infty = \|b_n\|_\infty$. Two forms of the completeness of $\mathcal{P}(G')$ can be used to complete the proof. We let $c_n = b_1 + \cdots + b_n$.

In the first proof, we notice that the sequence $S_{p'_1}, S_{p'_2}, \ldots$ of completely positive $C(\bar{X}')$-bimodule maps from $C^*(G', \bar{X}')$ to $M^*(R', X')$ is summable. The sum is also a completely positive $C(\bar{X}')$-bimodule map, so by Theorem 4.6 it is of the form $S_{p'}$ for a $p' \in \mathcal{P}(G')$. Then the function $b$ defined on $G$ by $b = p'(\cdot, (1,2))$ is in $\mathcal{B}(G)$ by Theorem 5.3. We also get $\|S_{p'_n - p'}\|_{cb} \geq \|S_{c_n - b}\|_{cb} \geq \|c_n - b\|_\infty$ by Theorem 5.3 and Theorem 5.2, so $\|c_n - b\|_\infty \to 0$. We need to prove that $\|c_n - b\| \to 0$ as $n \to \infty$.

To do this begin with $f \geq 0$ in $\mathcal{M}_c(G)$. Then Lemma 1.3 says that
$$\|\omega((c_n - b)f)\| \leq \|c_n - b\|_\infty \|\omega(f)\|.$$
Hence $T_{c_n}(\omega(f)) \to T_b(\omega(f))$ in $M^*(G)$. The $f$'s span a dense set in $M^*(G)$, and the $T_{c_n}$'s are uniformly bounded, so it follows that $T_{c_n} \to T_b$ pointwise on $M^*(G)$. Now the fact that the completely bounded operators on $M^*(G)$ implies that the sequence $T_{c_n}$ has a limit, $T'$ in the completely bounded sense, which is automatically also a pointwise limit on $M^*(G)$. Hence $T' = T_b$, so that $\|T_{c_n - b}\|_{cb} \to 0$, and by Theorem 5.1 that is equivalent to saying $\|c_n - b\| \to 0$ as $n \to \infty$.

For the other proof of completeness, we notice that $p'_1, p'_2, \ldots$ is summable in the $\mathcal{Q}$-essential supremum norm as functions on $G'$. Hence there is a Borel function $p'$ that is the sum in that norm. By the Dominated Convergence Theorem, $p' \in \mathcal{P}(G')$. Again we take $b = p'(\cdot, (1,2))$. Theorems 5.3 and 5.2 once again show that $\|c_n - b\|_\infty \to 0$, and we complete the proof as before.

Since $\mathcal{B}(G)$ is a Banach algebra, any closed subalgebra of it is a Banach algebra. Convergence in the completely bounded norm implies convergence in $L^\infty(\lambda^{\mathcal{Q}})$, so certain subalgebras are easily seen to be closed. Among these are $B(G)$, defined to be $\{b \in \mathcal{B}(G) : b$ is continuous$\}$, and $\mathcal{B}(G, X)$, defined to be the set of elements $b \in \mathcal{B}(G)$ such that $b|X$ is continuous and vanishes at $\infty$. The subalgebra $B(G, X)$ is defined to be $B(G) \cap \mathcal{B}(G, X)$.

**Theorem 6.2.** $B(G)$, $\mathcal{B}(G, X)$, and $B(G, X)$ are closed subalgebras of $\mathcal{B}(G)$ and hence Banach algebras.

## Section 7. Counter-examples.

The first example is a groupoid on which the linear span of the continuous positive definite functions is not complete and there exist continuous elements of $\mathcal{B}(G)$ that cannot be expressed as a difference of continuous positive definite functions.

32Let $X = \{(x,y) : (x,y)$ has polar coordinates $(r,\theta)$ with $0 \leq r \leq 1$, $\theta \in \{0, 1, 1/2, 1/3, \ldots\}\}\}$ and set $G = X \times \mathbb{Z}$. This is a bundle of groups, and $(x,n) + (x',n')$ is defined iff $x = x'$, and then it equals $(x, n+n')$. Write $\mathcal{P}(G)$ for the set of Borel positive definite functions on $G$ and $P(G)$ for the set of continuous elements of $\mathcal{P}(G)$. Let $\mathcal{B}(G)$ be the linear span of $\mathcal{P}(G)$, let $B_1(G)$ be the linear span of $P(G)$ and let $B(G)$ be the set of continuous elements of $\mathcal{B}(G)$. A bounded function $p$ is in $\mathcal{P}(G)$ iff it is a Borel function and $p(r,\theta,\cdot)$ is positive definite on $\mathbb{Z}$ for each point of $X$. Since positive definite functions on $\mathbb{Z}$ are in one-one correspondence with positive measures on $\mathbb{T}$ via the Fourier transform, we can also think of $\mathcal{P}(G)$ as consisting of Borel functions from $X$ to the positive measures on $\mathbb{T}$.

Define
$$p(r,\theta,n) = \begin{cases} e^{i\theta(1+r)n} & \text{if } r > 0 \\ 0 & \text{if } r = 0 \end{cases}$$
and
$$q(r,\theta,n) = \begin{cases} e^{i\theta(1-r)n} & \text{if } r > 0 \\ 0 & \text{if } r = 0 \end{cases}.$$

We can also think of these as taking values that are point masses at $e^{i\theta(1+r)}$ and $e^{i\theta(1-r)}$, or the 0 measure at the origin. We have $p - q \in B(G)$. Suppose that $u \in P(G)$ and $-u \leq p - q \leq u$ where the inequalities indicate the pointwise order in the space of measure-valued functions. This is the same as the natural order in $\mathcal{B}(G)$ in which elements of $\mathcal{P}(G)$ are positive. Since $p(r,\theta,\cdot)$ is the point mass at $e^{i\theta(1+r)}$, $u(r,\theta,\cdot)$ dominates the point mass at that point. By continuity, $u(0,0,\cdot)$ dominates the point mass at $e^{i\theta}$. This means that $u(0,0,\cdot)$ has infinite norm, so there is no such $u$. Thus we have a continuous element of $\mathcal{B}(G)$ that is not a difference of continuous positive definite functions.

With more effort, a worse example can be made. Choose $n$ angles, and begin with $p$ and $q$ restricted to the radii with those angles. The limit at the origin of both of them exists, the limits are the same, and it is a sum of $n$ point masses. To make elements of $P(G)$ we take that value at the origin and at all other points of $X$. Let $b$ be the difference of these elements of $P(G)$. Any element of $P(G)$ that dominates $b$ must have a value at the origin that dominates that sum of $n$ point masses. Observe that $b$ is 0 except on the original chosen radii, and that the total variation norm of each value of $b$ is at most 2.

Now partition the angles in $X$ into sets with $2^k$ elements, for $k = 1, 2, \ldots$, and use the construction just described to make elements $b_k$ in $B_1(G)$. Then let $b = \sum_{k=1}^{\infty} 2^{-k} b_k$. This converges in the completely bounded norm since each $b_k$ has completely bounded norm 2. Hence it also converges in uniform norm, so that $b \in B(G)$. Also $b$ is in the closure of $B_1(G)$. However, the domination arguments used above show that $b$ is not in $B_1(G)$.

The next example shows that locally compact groupoids can have unitary representations that are Borel but not continuous.

Consider an action of the integers on the circle by an irrational rotation, and form the transformation group groupoid, $G = \mathbb{T} \times \mathbb{Z}$. If $u$ is a unitary valued Borel function on $\mathbb{T}$, there is a unitary representation $U$ such that for all $\tau \in \mathbb{T}$, $u(\tau) = U(\tau, 1)$. If $u$ is not continuous, neither is $U$.

## References

[B–C–M–R]    L. W. Baggett, A. Carey, W. Moran, and A. Ramsay, Non-monomial multiplier representations of abelian groups, *J. Funct. Anal.* **97** (1990), 361–372.

[B–M–R]    L. Baggett, W. Mitchell, and A. Ramsay, Representations of the discrete Heisenberg group and cocycles of an irrational rotation, *Mich. J. Math*, **31** (1984), 263–273.

[B–R]    L. W. Baggett and A. Ramsay, A functional analytic proof of a selection lemma, *Canadian J. Math.*, **32** (1980), 441–448.